\documentclass[10pt]{article}
\usepackage[T1]{fontenc} 
\usepackage[utf8]{inputenc}

\usepackage{lmodern}

\usepackage{calrsfs}
\usepackage{amssymb}
\usepackage[cmex10]{amsmath}
\usepackage{amsthm}
\usepackage{latexsym}
\usepackage{epsfig}
\usepackage{graphicx}

\def\limsup{\mathop{\rm lim\,sup}}

\newcommand{\Prob}{\ensuremath{\mathsf{P}}}
\newcommand{\Ex}{\ensuremath{\mathsf{E}}}

\newtheorem{theorem}{Theorem}
\newtheorem{assumption}{Assumption}

\newtheorem{proposition}{Proposition}

\newtheorem{remark}{Remark}

\usepackage{appendix}

\usepackage[top=1in, bottom=1.25in, left=0.9in, right=0.9in]{geometry}

\usepackage[sort&compress]{natbib} 
\setcitestyle{authoryear,open={(},close={)}}

\usepackage[colorlinks=true,breaklinks=true,bookmarks=true,urlcolor=blue,
     citecolor=blue,linkcolor=blue,bookmarksopen=false,draft=false]{hyperref}

\def\EMAIL#1{\href{mailto:#1}{#1}}
\def\URL#1{\href{#1}{#1}}         

\begin{document}

\title{Admission and routing of soft real-time jobs to multiclusters: Design and comparison of index policies}


\author{Jos\'e Ni\~no-Mora 
\\ Department of Statistics \\
    Carlos III University of Madrid \\
     28903 Getafe (Madrid), Spain \\  \EMAIL{jnimora@alum.mit.edu}, \URL{http://orcid.org/0000-0002-2172-3983} }
\date{Published in \textit{Computers \& Operations Research}, vol.\ 39, 3431--3444, 2012 \\
DOI: \href{https://doi.org/10.1016/j.cor.2012.05.004}{10.1016/j.cor.2012.05.004}
}

\maketitle

\begin{abstract}%
Motivated by time-sensitive e-service applications, we consider the design of effective policies in a
Markovian model for the dynamic control of both admission and routing of a single class of real-time
transactions to multiple heterogeneous clusters of web servers, each having its own queue and server
pool. Transactions come with response-time deadlines, staying until completion if the latter are missed.
Per job rejection and deadline-miss penalties are incurred. Since computing an optimal policy is
intractable, we aim to design near optimal heuristic policies that are tractable for large-scale systems.
Four policies are developed: the static optimal Bernoulli-splitting (BS) policy, and three index policies,
based respectively on individually optimal (IO) actions, one-step policy improvement (PI), and restless
bandit (RB) indexation. A computational study demonstrates that PI is the best of such policies, being
consistently near optimal. In the pure-routing case, both the PI and RB policies are nearly optimal.
\end{abstract}%


\newpage
\section{Introduction} 
\label{s:intro}
Motivated by time-sensitive online data services (e-services) applications, this paper considers the design of effective policies for the dynamic  control of both admission and routing of a single class of real-time transactions to multiple clusters of web servers,
in the setting of a queueing model of Markov decision process (MDP) type. Service resources comprise $n$ parallel clusters of servers, where each  cluster $k \in \mathbb{K} \triangleq \{1, \ldots, n\}$ has its own dedicated 
queue with unlimited buffer space for holding transaction requests waiting or being processed, and a pool of $m_k$ identical servers, each working at a speed of $\mu_k$ (transactions per time unit). The system is heterogeneous in that both server speeds and server-pool sizes may differ across clusters.
Transactions arrive as a Poisson process with rate $\lambda$, with their execution times being independent and exponentially distributed with unit rate.
To ensure stability, it is assumed that the incoming load $\rho \triangleq \lambda / \sum_{k} m_k \mu_k$ is less than unity. 
Upon  arrival of a transaction request, the system controller immediately and irrevocably decides  whether to admit  or reject it. If admitted, the controller  further decides at once to which  cluster to route it, basing decisions on the history of previous queue lengths and actions.
Once in a cluster queue, pending transactions are scheduled  on a first-come first-served (FCFS) basis. 

Transactions are delay-sensitive of \emph{soft real-time} type, coming with timeliness requirements on their individual response or waiting times given in the form of \emph{soft deadlines} from their times of arrival, i.e., deadlines that can be missed without totally destroying the value of the transactions. 
Hence, transactions with missed deadlines remain in the system until completed.

The prime motivation for the study of such a system model comes from certain time-sensitive e-service applications, where it is desirable but not mandatory to meet given transaction deadlines, and which increasingly rely on   heterogeneous multicluster web-server architectures. Examples of such soft real-time applications include
e-commerce (see \cite{bertinietal10}), real-time databases, e.g., for online stock trading (see \cite{kaogarcmol96} and \cite{kangetal07}), and scientific computing (see \cite{zhu01} and \cite{plankenetal10}). See also \cite{heetal06}, which considers a similar model to the present one, yet focusing on the \emph{pure-routing} case, i.e., without admission control.

It should be mentioned that in other important applications, notably in call centers, transactions are subject to \emph{firm deadlines}, meaning that customers abandon if they are put on hold for too long. Use of index policies based on restless bandits for the resulting models with abandonment, or reneging, has been proposed in \cite{nmnetcoop07}. A thorough study of such models is the subject of another paper by the author, as yet unpublished.

The incorporation of admission control in the present model is motivated by the observation that, at times when the system is heavily congested,  admitting newly arriving transactions will almost guarantee that their deadlines will be missed. Hence, in such circumstances it may be beneficial to reject new transactions if doing so is not too costly.

\begin{figure}[!htb]
\centering
\includegraphics[height=2.5in]{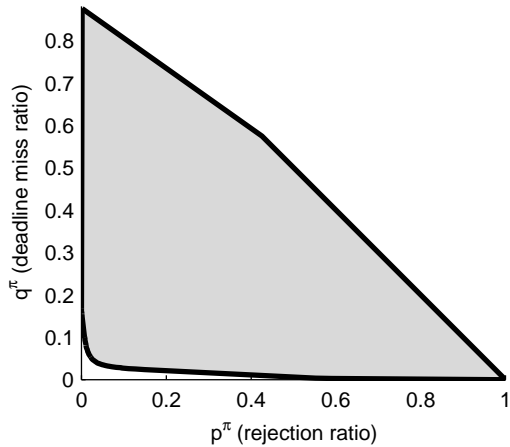}
\caption{The achievable performance region of $(p^\pi, q^\pi)$ pairs.}
\label{fig:achievReg}
\end{figure}

The main performance metrics to evaluate the Quality-of-Service (QoS) provided to users under a given control policy $\pi$ in such a system are the \emph{rejection ratio} $p^\pi$, which is the long-run average fraction of transactions that are rejected upfront, and the \emph{deadline miss ratio} $q^\pi$, which is the long-run average fraction of transactions that are admitted only to miss their deadlines. While one would like to have both $p^\pi$ and $q^\pi$ small, it is intuitively clear that there is trade-off between them. We can visualize such a trade-off as the lower boundary in the \emph{region of achievable performance} of  $(p^\pi, q^\pi)$ pairs, spanned under all admissible policies $\pi$, which is displayed in Figure \ref{fig:achievReg} for the instance with $n = 2$ clusters with server-pool sizes $\mathbf{m} = (4, 8)$, server speeds $\boldsymbol{\mu} = (5, 3)$, and $\rho = 0.95$. 

A common approach to policy design in such systems is to specify in advance upper bounds $\bar{p}$ and $\bar{q}$ on the metrics $p^\pi$ and $q^\pi$, and then try to  construct a policy $\pi$ that satisfies the resulting constraints $p^\pi \leqslant \bar{p}$ and $q^\pi\leqslant \bar{q}$. Another common approach is to aim to minimize one of the metrics, say $p^\pi$, subject to an upper bound constraint on the other, e.g., $q^\pi\leqslant \bar{q}$.

In contrast, in this work we aim to jointly optimize the pair $(p^\pi, q^\pi)$, by setting economic incentives to the service provider, considering a \emph{performance objective} to minimize of the form 
$R p^\pi + C q^\pi$, where $R > 0$ represents a cost per rejected transaction, and $C > 0$ is a cost per missed deadline. Notice that the  pure-routing case considered in \cite{heetal06} can be seen as the special case of the present model where the rejection cost is $R$ is sufficiently large.

Since minimizing the cost $R p^\pi + C q^\pi$ is equivalent to maximizing the net reward $R (1~-~p^\pi) - C q^\pi$, such an objective models a Service Level Agreement (SLA) contract, which sets the fee $R$ to be collected by the service provider from each admitted user, and incorporates an economic compensation scheme by which a refund $C$ is issued to a user whose QoS target deadline is missed.
While the above \emph{performance optimization problem} can be formulated as a countable-state average-cost MDP, computing an optimal policy by solving numerically the corresponding \emph{dynamic programming} (DP) equations is computationally intractable, even for a moderate number of clusters. 
Hence, the main aim of this paper is to design high-performance heuristic policies that can be computed with low complexity in large-scale systems.

In the sequel we will focus, with no loss of generality, on the case $C = 1$, which results by 
normalizing the above  objective by dividing it by $C$.  Further, we will use interchangeably the terms ``queue'' and ``cluster,'' since each cluster has its own queue, and we will use the term ``job'' instead of ``transaction.''

\subsection{Related prior work}
\label{s:rpw}
Most related prior  work on soft-real time systems has focused on pure-routing models, and has only
considered \emph{static} policies, which do not use system state information. Static policies decompose routing decisions into a \emph{load balancing} or \emph{workload allocation} rule, which determines the proportion of jobs to be assigned to each queue, and a \emph{job dispatching} rule, which specifies where to send each job. An admission policy may be incorporated by fixing the proportion of jobs to be rejected.
A simple static policy is a \emph{Bernoulli splitting} (BS),  which splits the arrival stream according to fixed probabilities.
\cite{kallmescass95} addresses computation of the optimal BS in the special case of the present model, including admission and routing, where each queue has a single server and the rejection and deadline miss penalties are equal ($R = C = 1$).
This calls for solving a convex separable \emph{nonlinear program} (NLP). 
\cite{heetal06} shows how to compute the optimal workload allocation for the special pure-routing case of the present model, considering BS as well as round-robin job dispatching.

As for dynamic policies, which use some degree of state information, prior work has addressed the \emph{performance evaluation} of the deadline-miss ratio for some ad hoc routing policies, in homogeneous systems with parallel  single-server queues. 
\cite{kaogarcmol96} considers routing policies that use deadline information, which are shown to outperform symmetric BS routing.
\cite{zhu01} reports on a comparative simulation study of four routing policies for the same model, one of which is the \emph{shortest queue routing} (SQR) rule. It is shown that SQR is the best policy among those considered.

In different models from the one considered here, index policies have been proposed for routing to parallel queues based on one-step \emph{policy improvement} (PI) over the optimal BS (see \cite{krish90}), and also based on \emph{restless bandit} (RB) indexation, introduced in general in \cite{whit88} and first applied in \cite{nmmp02} to problems of admission control and routing to parallel queues (see also \cite{nmnetcoop07}).

Effective heuristic index policies have also been proposed and investigated for dynamic routing in certain call-center models that incorporate both multiple customer classes and multiple server pools with different skills. See, e.g., \cite{bhulai09}, which derives a PI-based policy that is shown to perform extremely well, and \cite{gurvWhitt10}, where the proposed policies can achieved desired QoS levels.
Note, however, that, even in the case of a single customer class the models considered in such papers are fundamentally different from the one addressed here, not only because they do not incorporate rejections and consider either no deadlines or firm deadlines (with abandonment) rather than soft deadlines, but more importantly because they assume that customers (from each class) join a \emph{common queue} upon arrival, to be \emph{later} routed to available servers. In contrast, in the present model customer transactions that are admitted are immediately routed upon arrival to one of \emph{multiple parallel queues}, each corresponding to a multiserver cluster.

\subsection{Extending prior work}
\label{s:epw}
This paper extends the aforementioned research in the following ways:  (1) the  model incorporates dynamic control of admission, rather than considering only static or only pure-routing policies; 
(2) instead of addressing the performance evaluation of given ad hoc
policies, the aim is \emph{performance optimization}, i.e., to \emph{design} heuristic policies that can be computed with low complexity and come close to optimizing the performance objective of concern; (3) the focus is on systems with parallel \emph{heterogeneous} multi-server queues, motivated by the prevalence of heterogeneity in real multicluster systems; (4) the computation of the optimal BS is extended from the special case $R = 1$ (recall that $C \equiv 1$) to the general case where $R$ may differ from $1$; (5) new dynamic index policies  type are introduced, based on three different methods (individually optimal decisions, one-step policy improvement, and restless bandits); and, finally, all the heuristic policies considered are benchmarked against the optimal performance in a numerical study.

We will focus on the intuitively appealing class of \emph{index policies}, which determine admission and routing actions on the basis of $n$ numeric \emph{routing indices}, where index 
 $\nu_k(x_k)$ is attached to queue $k \in \mathbb{K}$ as a
nondecreasing
function of its current \emph{state} $x_k$, the latter giving the number of jobs waiting or being served in the cluster's queue.
The  resulting routing policy  sends an admitted job to a  queue with currently
 lowest index value, breaking ties arbitrarily.

Such routing policies aim to approximate the optimal routing policy, which is shown below in Section \ref{s:toari} to have the following form: a newly arriving job that finds the system in the joint state $\mathbf{x} = (x_k)$, if admitted, is assigned 
to a queue $k$ attaining the minimum value of a certain quantity $\nu_k^*(\mathbf{x})$, which is thus an \emph{optimal routing index} which depends on the joint state $\mathbf{x}$, rather than on the state $x_k$ of queue $k$ alone.
The optimal index $\nu_k^*(\mathbf{x})$ measures the long-run expected increment in total costs accrued that results by routing arrivals in state $\mathbf{x}$ to queue $k$, given that all other actions are taken optimally.

It is further shown in Section \ref{s:toari} that it is optimal to reject an arrival finding the system in the joint state $\mathbf{x} = (x_k)$ iff $R \leqslant \gamma^*(\mathbf{x})$, where 
$\gamma^*(\mathbf{x})$ is the \emph{optimal global admission index} defined by $\gamma^*(\mathbf{x}) \triangleq \min_k \nu_k^*(\mathbf{x})$.

Hence, for consistency with such optimal admission policies, given heuristic routing indices $\nu_k(x_k)$ we will base admission decisions on a corresponding
global admission index $\gamma(\mathbf{x})  \triangleq \min_k \nu_k(x_k)$.
 The resulting policy admits an arrival finding the system  in
 state $\mathbf{x}$ if $\gamma(\mathbf{x}) < R$, and rejects it otherwise. 

This raises the issue of how to define such routing indices so that they can be computed quickly and yield high-performance policies.
An insight on which we will draw is that the routing index $\nu_k(x_k)$ should be a proxy measure for the optimal index $\nu_k^*(\mathbf{x})$, which, as pointed out above, is the incremental cost of routing an arrival to queue $k$ in state $\mathbf{x}$. 

To design the routing indices $\nu_k(x_k)$ 
we will deploy three alternative general methods  that have been proposed in the literature for obtaining priority indices. Two such methods
correspond to different assessments of the  incremental cost of assigning an arrival to a queue.
At one extreme, the \emph{individually optimal} (IO) method \emph{underestimates} such a cost by disregarding
future arrivals. At the opposite extreme, the RB method \emph{overestimates} such a cost by taking into account all future arrivals, as 
though there were no other queues to share the load.  
Such methods yield  routing indices $\nu_k^{\textup{IO}}(x_k)$ and $\nu_k^{\textup{RB}}(x_k)$ that can be evaluated separately for each cluster $k$, as they depend only on its individual  parameters.
The 
third method we will consider, which is to carry out one PI  step of the policy iteration algorithm for MDPs, starting from
the optimal BS, gives routing indices $\nu_k^{\textup{PI}}(x_k)$ that depend on the parameters of all clusters.
Such a method has been used in \cite{krish90} to obtain a heuristic index policy for job assignment to parallel multiserver queues with a minimum average response-time objective.

To the best of the author's knowledge, none of the above policy-design methods has been applied in prior work to the present model. Only the optimal BS method has been applied to special cases: in \cite{kallmescass95} for the  case of single-server queues $m_k \equiv 1$, and in \cite{heetal06} for the pure-routing case.

The main objective of the present work is to deploy the BS, PI, IO, and RB policy design methods in the present model, addressing the efficient computation of the resulting policies, and to identify the best among such policies through a numerical study.

The remainder of the paper is organized as follows. 
Section \ref{s:cao} describes the model and formulates the optimal control problem. 
Section \ref{s:model1io} addresses the efficient evaluation of the IO routing index in the cases of constant, uniform, and exponentially distributed job deadlines.
Section \ref{s:model2bs} discusses the  PI routing index, as well as the prerequisite computation of the optimal BS.
Section \ref{s:mod2mp} considers the RB routing index.
Section \ref{s:ne} reports on the results of a numerical study benchmarking the above policies, against each other and against the optimal performance.
 Section \ref{s:concl} concludes.

\section{The model}
\label{s:cao}
\subsection{Model description and problem formulation}
\label{s:md}
Consider
 a  system with $n$ parallel clusters of servers, with 
cluster $k \in \mathbb{K}$ 
having its own queue with unlimited buffer space and $m_k$ identical exponential servers, each working at  rate $\mu_k$.
Jobs arrive as
a Poisson stream with 
rate $\lambda$, with job interarrival and service times being independent.
Upon  arrival of a job, it is decided  (i) whether to 
admit  or reject it;
 and,  if admitted, (ii) to which
queue to  route it. 

Jobs come with a (relative) deadline $\tau$, giving the time interval from arrival to deadline expiration, which is independent of arrival and service times, being 
drawn in i.i.d.\ fashion from a general cumulative distribution function (CDF) $F(t)$ with finite mean $\Ex[\tau] = 1/\theta$. 
We will consider the cases of deterministic, uniform and exponential deadlines, 
 focusing on the case of \emph{deadlines to the end of service}, where $\tau$ represents an individual \emph{response-time deadline}, although the ensuing analyses readily extend to the case of deadlines to the beginning of service.
Admitted jobs whose deadlines are missed stay nevertheless in the system until service completion.
A cost of $R > 0$ is incurred
per rejected job, whereas a cost of $C = 1$ is incurred per missed deadline.

Denote by $X_k(t)$ the state
 of queue $k$ at time $t$,  giving the number of jobs its queue holds waiting or being served, and by $a_k(t) \in \{0,
1\}$ the binary \emph{action}  that  is set to $0$  if
a job arriving at time $t$  \emph{would be routed to queue}
$k$. One may imagine that each queue $k$ has its own \emph{entry
  gate} and \emph{gatekeeper}, who 
 can open ($a_k(t) = 0$) or shut ($a_k(t) =
1$) the queue's entry gate. 
The requirement that a job can be routed to at most one queue is thus formulated by 
the constraint that at most one gate can be open at any time:
 \begin{equation}
 \label{eq:spconstr}
 \sum_{k \in \mathbb{K}} \big(1-a_k(t)\big)  \leqslant 1, \quad t \geqslant 0.
 \end{equation}

The following condition will be assumed to hold. 
\begin{assumption}
\label{ass:modelii1}
At each queue $k$, the queue discipline is FCFS.
\end{assumption}

Assumption \ref{ass:modelii1}, which excludes use of alternative queue disciplines such as \emph{earliest deadline first}, is needed to ensure that the  response-time distribution of a job assigned to a queue $k$ is determined only by the queue state
$x_k$ found on arrival. 
Denoting by $T_k(x_k)$ a random variable with such a  \emph{conditional response-time distribution}, and letting $\tau$ be the job's deadline, the expected deadline-miss cost incurred by a job assigned to queue $k$ in state $x_k$ is  
\begin{equation}
\label{eq:rkikmiiidef}
P_k^{\textup{miss}}(x_k) \triangleq \Prob\left\{T_k(x_k) > \tau\right\} =  \int_0^\infty \Prob\left\{T_k(x_k) > t\right\} \, dF(t).
\end{equation}

We will find convenient to formulate the lump costs incurred due to job rejections and deadline misses in terms of continually accruing costs, by considering 
that the total cost rate incurred per unit time 
 when the joint action $\mathbf{a} = (a_k)_{k \in \mathbb{K}}$ is taken in the joint state
$\mathbf{x} = (x_k)_{k \in \mathbb{K}}$  is given by 
\begin{equation}
\label{eq:obj1}
 \lambda R \Bigg[1 - \sum_{k \in \mathbb{K}} (1-a_k)\Bigg]  + \sum_{k \in \mathbb{K}}  \lambda P_k^{\textup{miss}}(x_k) (1-a_k).
\end{equation}
We will also find convenient to reformulate  (\ref{eq:obj1}) to separable form as
\begin{equation}
\label{eq:obj}
-(n-1)  \lambda R + \sum_{k \in \mathbb{K}} \big\{\lambda  R a_k + \lambda P_k^{\textup{miss}}(x_k) (1-a_k)\big\}.
\end{equation}

Consider the  
problem of 
finding an \emph{average-cost optimal policy}, which minimizes the long-run expected average 
cost per unit time. Disregarding the constant $-(n-1)  \lambda R$ in (\ref{eq:obj}), we can formulate such a problem as
\begin{equation}
\label{eq:pa}
\min_{\boldsymbol{\pi} \in \boldsymbol{\Pi}} \limsup_{T \to \infty} \frac{1}{T}
\Ex_{\mathbf{x}^0}^{\boldsymbol{\pi}}\Bigg[\int_0^T  \sum_{k \in \mathbb{K}} \Big\{ \lambda  R 
 a_k(t) + \lambda P_k^{\textup{miss}}\big(X_k(t)\big) \big(1-a_k(t)\big) \Big\} \, dt\Bigg],
\end{equation}
where ``$\limsup$'' is the limit superior, 
$\Ex_{\mathbf{x}^0}^{\boldsymbol{\pi}}[\cdot]$ denotes expectation under
policy $\boldsymbol{\pi}$ conditioned on the initial joint state being equal
to $\mathbf{x}^0 = (x_k^0)_{k \in \mathbb{K}}$, and $\boldsymbol{\Pi}$ is the class of \emph{admissible  policies}, which are allowed to depend only on the history of system states and actions encountered in the past, and may be randomized.

By standard results, problem (\ref{eq:pa}) has an optimal policy that is \emph{stationary deterministic}, i.e., in which the action to take at every state $\mathbf{x}$ is a deterministic function of $\mathbf{x}$. Both optimal policies and the optimal cost rate are independent of the initial state, and they are characterized by the \emph{DP average-cost optimality equation} (see (\ref{eq:dpeqjor}) below). 
Yet,  finding an optimal policy by numerical solution of such equations is, due to the curse of dimensionality, computationally intractable for instances with even a moderate number $n$ of queues. Hence, we aim to design and compare heuristic index policies that can be computed with low complexity and perform well, as discussed in Section \ref{s:intro}.

 \subsection{Charging  deadline-miss costs via action-independent holding costs}
 \label{s:rjosor}
Application of two of the index design methods considered in this paper, PI and RB (see Sections \ref{s:model2bs} and  \ref{s:mod2mp}), would be considerably simplified if, rather than 
charging  the deadline-miss costs incurred at each queue $k$ in terms of action-dependent cost rates as in formulation (\ref{eq:pa}), they were charged  in terms of action-independent holding costs, which continually accrued at a certain rate
 $h_k(x_k)$ per unit time when the queue holds $x_k$ jobs.

Hence, to prepare the ground for the analyses in Sections \ref{s:model2bs} and  \ref{s:mod2mp}, we next address the reformulation of problem (\ref{eq:pa}) to an equivalent problem where deadline-miss costs are accounted for via action-independent holding cost rates  $h_k(x_k)$.
For such a purpose, we start by formulating the \emph{DP optimality equations} characterizing the optimal policies for average-cost problem (\ref{eq:pa}), which are stated in terms of variables $v$ and $b(\mathbf{x})$: for every system state $\mathbf{x}$,
\begin{equation}
 \label{eq:dpeqjor}
 v = \min
 \begin{cases}
 \displaystyle n \lambda R  -  \sum_{l \in \mathbb{K}} \bar{\mu}_l(x_l) \Delta_l b(\mathbf{x}) \\
\displaystyle  (n-1) \lambda R + \lambda P_k^{\textup{miss}}(x_k) + \lambda \Delta_k b(\mathbf{x} + \mathbf{e}_k) -  \sum_{l \in \mathbb{K}\colon x_l \geqslant 1} \bar{\mu}_l(x_l) \Delta_l b(\mathbf{x}),  k \in \mathbb{K},
 \end{cases}
 \end{equation}
where we write $\bar{\mu}_k(x_k) \triangleq \min(x_k, m_k) \mu_k$
 and
 $\Delta_k b(\mathbf{x}) \triangleq b(\mathbf{x}) - b(\mathbf{x}-\mathbf{e}_k)$, $\mathbf{e}_k$ being the $k$th unit coordinate vector in $\mathbb{R}^n$.  By standard results, the DP equations (\ref{eq:dpeqjor}) have a solution 
 $v^*$ and $b^*(\mathbf{x})$, which is unique up to a constant additive term for the $b^*(\mathbf{x})$, where $v^*$ is the optimal average cost per unit time and $b^*(\mathbf{x})$ is an optimal relative cost function.

Let us now define the  holding cost rate $h_k(\cdot)$ for each queue $k \in \mathbb{K}$ by
\begin{equation}
\label{eq:ridef}
h_k(0) \triangleq 0, \quad h_k(x_k) \triangleq \bar{\mu}_k(x_k) P_k^{\textup{miss}}(x_k-1), \quad x_k \geqslant 1,
\end{equation}
and consider the following modified admission control and routing problem, where two types of costs are incurred: rejection costs, as before, and holding costs, which accrue at each queue $k$ at rate $h_k(x_k)$ while it holds $x_k$ jobs:
\begin{equation}
\label{eq:pa2}
\min_{\boldsymbol{\pi} \in \boldsymbol{\Pi}} \limsup_{T \to \infty} \frac{1}{T}
\Ex_{\mathbf{x}^0}^{\boldsymbol{\pi}}\Bigg[\int_0^T  \sum_{k \in \mathbb{K}} \Big\{ R  \lambda a_k(t) + h_k\big(X_k(t)\big)\Big\} \, dt\Bigg].
\end{equation}

The optimal policies for the modified problem (\ref{eq:pa2}) are characterized by its DP optimality equations: for every system state $\mathbf{x}$, 
\begin{equation}
 \label{eq:dpeqnsmodel}
\widetilde{v} = \min
 \begin{cases}
 \displaystyle n \lambda R + \sum_{l \in \mathbb{K}} h_l(x_l) -  \sum_{l \in \mathbb{K}} \bar{\mu}_l(x_l) \Delta_l \widetilde{b}(\mathbf{x}) \\
\displaystyle (n-1) \lambda R + \sum_{l \in \mathbb{K}} h_l(x_l)  + \lambda \Delta_k \widetilde{b}(\mathbf{x} + \mathbf{e}_k) -  \sum_{l \in \mathbb{K}\colon x_l \geqslant 1} \bar{\mu}_l(x_l) \Delta_l \widetilde{b}(\mathbf{x}), \, k \in \mathbb{K}.
 \end{cases}
 \end{equation}
Equations (\ref{eq:dpeqnsmodel}) have a solution $\widetilde{v}^*$ and $\widetilde{b}^*(\mathbf{x})$, which is unique up to a constant additive term for the $\widetilde{b}^*(\mathbf{x})$, where 
$\widetilde{v}^*$ is the optimal average cost per unit time and $\widetilde{b}^*(\mathbf{x})$  is an optimal relative cost function.

The following result ensures that, charging costs at each queue using the holding cost rates defined in (\ref{eq:ridef}), we obtain an equivalent problem.
\begin{proposition}
\label{pro:rjosor}
Problems $(\ref{eq:pa})$ and $(\ref{eq:pa2})$ are equivalent, having the same optimal policies and average cost per unit time.
\end{proposition}
\begin{proof}
For each queue $k$, define the function $\widetilde{C}_k(\cdot)$ by the recursion 
 \begin{equation}
 \label{eq:rkrecs}
\widetilde{C}_k(0) \triangleq 0, \quad \widetilde{C}_k(x_k+1) \triangleq \widetilde{C}_k(x_k) + P_k^{\textup{miss}}(x_k), \quad x_k \geqslant 0.
 \end{equation}
The result follows by checking that
$v^* \triangleq \widetilde{v}^*$  and  $b^*(\mathbf{x}) \triangleq \widetilde{b}^*(\mathbf{x}) - \sum_k \widetilde{C}_k(x_k)$ satisfy the
DP equations (\ref{eq:dpeqjor}) iff $\widetilde{v}^*$ and the $\widetilde{b}^*(\mathbf{x})$ satisfy  (\ref{eq:dpeqnsmodel}). 
\end{proof}

\subsection{The structure of optimal admission and routing policies}
\label{s:toari}
In light of (\ref{eq:dpeqnsmodel}) and Proposition \ref{pro:rjosor}, an optimal policy for problem (\ref{eq:pa}) can be formulated in terms of the relative cost increments $\Delta \widetilde{b}^*(\mathbf{x})$ for 
problem (\ref{eq:pa2}): reject an arrival in state $\mathbf{x}$ if $R \leqslant \min_k \Delta_k \widetilde{b}^*(\mathbf{x} + \mathbf{e}_k)$; otherwise, admit the arrival and route it to a queue $k$ with minimum $\Delta_k \widetilde{b}^*(\mathbf{x} + \mathbf{e}_k)$.
Hence, optimal admission actions are characterized by a global admission index, given by $\gamma^*(\mathbf{x}) \triangleq \min_k \Delta_k \widetilde{b}^*(\mathbf{x} + \mathbf{e}_k)$, whereas optimal routing actions are characterized by indices $\nu_k^*(\mathbf{x}) \triangleq \Delta_k \widetilde{b}^*(\mathbf{x} + \mathbf{e}_k)$ that depend on the system state $\mathbf{x}$.

Given the unavailability of such optimal admission and routing indices, this paper compares, in terms of the performance of the resulting policies, several alternative  ways of defining a global admission index $\gamma(\mathbf{x})$ and routing indices $\nu_k(x_k)$ that can be evaluated with low computational complexity. 

\section{IO index}
\label{s:model1io}
The simplest dynamic routing policy we consider is to route an arrival to a
queue where the probability of missing its deadline, conditioned on the number of jobs found on such a queue upon arrival,  is minimal. The resulting IO routing index is given by such a conditional deadline-miss probability, being $\nu_k^{\textup{IO}}(x_k) \triangleq P_k^{\textup{miss}}(x_k)$, where $P_k^{\textup{miss}}(x_k)$ is defined in (\ref{eq:rkikmiiidef}).
Note that $\nu_k^{\textup{IO}}(0) = \cdots = \nu_k^{\textup{IO}}(m_k~-~1)$, and that $\nu_k^{\textup{IO}}(x_k)$  converges monotonically up to unity as $x_k \to \infty$:
\begin{equation}
\label{eq:limiio}
\bar{\nu}_k^{\textup{IO}}(\infty) = \lim_{x_k \to \infty} \nu_k^{\textup{IO}}(x_k) = \lim_{x_k \to \infty}  P_k^{\textup{miss}}(x_k) = 1.
\end{equation}
Hence, the admission policy based on the IO routing index, having admission index $\gamma^{\textup{IO}}(\mathbf{x}) = \min_k \nu_k^{\textup{IO}}(x_k)$, will not reject any job in the case $R \geqslant 1$.

We next address the efficient evaluation of $P_k^{\textup{miss}}(x_k)$ in the cases of constant, uniform, and exponential response-time deadlines. Since the focus is on a single queue, the  label $k$ is dropped below from the notation, writing, e.g., $P^{\textup{miss}}(x)$. 

\subsection{Constant response-time deadlines}
\label{s:model1iodd}
In the case where jobs have a constant response-time deadline $\tau \equiv t$, $P^{\textup{miss}}(x)$ equals
$P_m(x; t) \triangleq \Prob\{T(x) > t\}$,  where $m$ is the number of servers in the queue and $T(x)$ denotes a random variable distributed as the response time (under FCFS) of an arrival who finds $x$ jobs present in the queue. To evaluate $P_m(x; t)$, 
we will need the following probabilities, where $A_m(j) \sim \mathrm{Erlang}(j, m \mu)$:
\begin{equation}
\label{eq:qnjdef}
Q_m(j; t) \triangleq \Prob\{A_m(j) > t\} = e^{-m \mu t} \sum_{l=0}^{j-1} \frac{(m \mu t)^l}{l!}, \quad j = 1, 2, \ldots
\end{equation}
For fixed $m$ and $t$,
the $Q_m(j; t)$ are efficiently evaluated by the following  second-order linear recursion, where
$\Delta Q_m(j; t) \triangleq  Q_m(j; t) -  Q_m(j-1; t)$:
\begin{equation}
\label{eq:pmjrec}
\begin{split}
Q_m(1; t) & = e^{-m \mu t}, \quad \Delta Q_m(2; t) =   m \mu t Q_m(1; t) \\
  \Delta Q_m(j+1; t) & = (m \mu t /{j}) \Delta Q_m(j; t), \quad j \geqslant 2.
\end{split}
\end{equation}

The following result evaluates $P_m(x; t)$ in terms of $Q_m(\cdot, t)$ and $Q_{m-1}(\cdot, t)$. Along with (\ref{eq:pmjrec}), it allows us to compute $P^{\textup{miss}}(0)$, \ldots, $P^{\textup{miss}}(x)$ in $O(x)$ time.

\begin{proposition}
\label{pro:probtmifor}
\begin{itemize}
\item[\textup{(a)}] $P_1(x; t) = Q_1(x+1; t)$, for $x \geqslant 0$.
\item[\textup{(b)}] If $m \geqslant 2$, then
$P_m(x; t)  = e^{-\mu t}$ for $0 \leqslant x < m$; and, for $x \geqslant m$,
\[
P_m(x; t) = 
Q_m(x-m+1; t) + e^{-\mu t}  \frac{1 - Q_{m-1}(x-m+1; t)}{(1-1/m)^{x-m+1}}.
\]
\end{itemize}
\end{proposition}
\begin{proof}
(a) If $m = 1$,  $T(x) \sim \mathrm{Erlang}(x~+~1,~\mu)$, which yields the result.

(b) If $m \geqslant 2$, then
$P_m(x; t)  = e^{-\mu t}$ for $0 \leqslant x < m$, as in this case the arriving job finds some free server. Consider now the case $x \geqslant m$. Then, the conditional response time is decomposed as
$T(x) \sim A + \xi$, where $A \sim \mathrm{Erlang}(x-m+1, m \mu)$ and  $\xi \sim \mathrm{Exp}(\mu)$ are independent random variables giving the job's conditional waiting time and its service time, respectively.
Hence, denoting by $f_A(t)$ the density of $A$, and writing $j = x-m+1$,
 \begin{align*}
 P_m(x; t) & = \Prob\{A +  \xi > t\} \\
& =  \Prob\{A +  \xi > t \, | \, A > t\} \Prob\{A > t\}  + \int_0^t \Prob\{A +  \xi > t \, | \, A = s\} f_A(s) \, ds \\
 & = \Prob\{A > t\} + \int_0^t \Prob\{\xi > t - s\} f_A(s) \, ds  = Q_m(j; t) + \int_0^t e^{-\mu (t -s)} f_A(s) \, ds.
 \end{align*}
 
 On the other hand, letting $A' \sim \mathrm{Erlang}\big(j, (m-1) \mu\big)$, we have
 \begin{align*}
 \int_0^t e^{-\mu (t - s)} f_A(s) \, ds & = \int_0^t e^{-\mu (t - s)} (m \mu)^j \frac{s^{j-1}}{(j-1)!} e^{-m \mu s} \, ds  \\
& = \frac{e^{-\mu t}  }{(1-1/m)^{j}} \int_0^t  \big((m-1) \mu\big)^j \frac{s^{j-1}}{(j-1)!} e^{-(m-1) \mu s} \, ds \\
 & = \frac{e^{-\mu t}  }{(1-1/m)^{j}} \int_0^t f_{A'}(s) \, ds  = e^{-\mu t}  \frac{1 - Q_{m-1}(j; t)}{(1-1/m)^{j}}.
\end{align*}
\end{proof}

\subsection{Uniform response-time deadlines}
\label{s:ious}
Assuming that jobs have response-time deadlines $\tau \sim \mathrm{Unif}[t_{1}, t_{2}]$, $P^{\textup{miss}}(x) = P_m^*(x)/\big(t_{2} - t_{1})$, where 
$P_m^*(x) \triangleq
\int_{t_{1}}^{t_{2}} P_m(x; t) \, dt$. To evaluate $P_m^*(x)$, we will use quantities   
$Q_m^*(j) \triangleq
\int_{t_{1}}^{t_{2}} Q_m(j; t) \, dt$
and $R_{m-1}^*(j) \triangleq
\int_{t_{1}}^{t_{2}} e^{-\mu t} Q_{m-1}(j; t) \, dt$.

Both $Q_m^*(j)$ and $R_{m-1}^*(j)$ are readily evaluated in terms of the $Q_m(\cdot; \cdot)$ and $Q_{m-1}(\cdot; \cdot)$ using
integration by parts, which yields
\begin{equation}
\label{eq:qnj}
Q_m^*(j) = t_{2} Q_m(j; t_{2}) -
t_{1} Q_m(j; t_{1}) + \frac{j}{m
  \mu}\big[Q_m(j+1;  t_{1}) - Q_m(j+1;  t_{2})\big],
\end{equation}
and
\begin{equation}
\label{eq:rnj}
R_{m-1}^*(j) = \frac{1}{\mu}\left[e^{-\mu t_{1}} Q_{m-1}(j;
  t_{1}) - e^{-\mu t_{2}} Q_{m-1}(j;
  t_{2}) - \frac{Q_{m}(j; t_{1}) -
    Q_{m}(j; t_{2})}{(1 + 1/(m-1))^j}\right].
\end{equation}

The next result, which follows from Proposition \ref{pro:probtmifor}, evaluates $P_m^*(x)$ in terms of $Q_m^*(\cdot)$ and $R_{m-1}^*(\cdot)$.
Along with (\ref{eq:qnj}), (\ref{eq:rnj}), and the results in Section \ref{s:model1iodd}, this  allows us to evaluate $P^{\textup{miss}}(x)$, for $i = 0, \ldots, j$, in  $O(j)$ time.

\begin{proposition}
\label{pro:probunifs}
\begin{itemize}
\item[\textup{(a)}] $P_1^*(x) = Q_1^*(x+1)$, for $x \geqslant 0$.
\item[\textup{(b)}] If $m \geqslant 2$, then 
$P_m^*(x)  = \big(e^{-\mu t_{1}}-e^{-\mu t_{2}}\big)/\mu$ for $0 \leqslant x < m$; for $x \geqslant m$,
\[
P_m^*(x) = 
Q_m^*(x-m+1) + \frac{P_m^*(0) - R_{m-1}^*(x-m+1)}{(1-1/m)^{x-m+1}} .
\]
\end{itemize}
\end{proposition}

\subsection{Exponential  response-time deadlines}
\label{s:model1ioed}
Assuming that jobs have response-time deadlines $\tau \sim \mathrm{Exp}(\theta)$, $P^{\textup{miss}}(x) = \theta P_m^*(x; \theta)$, where
$P_m^*(x; s) \triangleq \int_0^\infty e^{-s t} P_m(x; t) \, dt$. To evaluate $P_m^*(x; \theta)$ we will use quantities 
$Q_m^*(j; s) \triangleq \int_0^\infty e^{-s t} Q_m(j; t) \, dt$, which, by standard results, are evaluated in closed form as
\[
Q^*_m(j; s) = \frac{1}{s} \Bigg[1 - \bigg(\frac{m \mu}{s + m \mu}\bigg)^j\Bigg], \quad j = 1, 2, \ldots
\]

The next result follows from Proposition \ref{pro:probtmifor} by taking Laplace transforms.
\begin{proposition}
\label{pro:probtmifor2}
\begin{itemize}
\item[\textup{(a)}] $P_1^*(x; \theta) = Q_1^*(x+1; \theta)$, for $x \geqslant 0$.
\item[\textup{(b)}] If $m \geqslant 2$, then
$P_m^*(x; \theta)  = 1/(\mu + \theta)$ for $0 \leqslant x < m$; and, for $x \geqslant m$,
\[
P_m^*(x; \theta) = 
Q_m^*( x-m+1; \theta) + \frac{\displaystyle \frac{1}{\mu + \theta} - Q_{m-1}^*( x-m+1; \mu + \theta)}{(1-1/m)^{ x-m+1}}.
\]
\end{itemize}
\end{proposition}

\section{PI index}
\label{s:model2bs}
We next consider the PI method (see  \cite{krish90}), which has not been applied before to the present model. 
This method consists of two stages: (1) computing the optimal BS of the arrival stream; and (2) carrying out one step of the PI algorithm for MDPs, starting from the optimal BS.

\subsection{First stage of the PI method: computing the optimal BS}
\label{s:fspiim}
The first stage is to compute the optimal BS of the arrival stream. In a BS,  the actions taken upon job arrivals are independently drawn according to fixed probabilities: a job is rejected with probability $\lambda_0/\lambda$, and is assigned to queue $k$ with probability $\lambda_k/\lambda$, for $k \in \mathbb{K}$. 
Hence, the input to each queue $k$ is a Poisson process with rate $\lambda_k$, with the queues evolving independently. 
Consider the resulting M/M/$m_k$ queue $k$, which has offered load $r_k(\lambda_k) \triangleq 
\lambda_k/\mu_k$ and utilization factor $\rho_k(\lambda_k) \triangleq 
\lambda_k/(m_k \mu_k)$. Assuming that $\rho_k(\lambda_k) < 1$, 
let $\widetilde{T}_k(\lambda_k)$ be a random variable with the queue's steady-state sojourn-time distribution, whose CDF $S_k(t; \lambda_k) \triangleq \Prob\{\widetilde{T}_k(\lambda_k) \leqslant t\}$ is given by (see  \citet[p.\ 72]{grossetal08})
\begin{equation}
\label{eq:wktlambdak}
S_k(t; \lambda_k)  = 
1 - e^{-\mu_k t} - G_k(\lambda_k)  \big(e^{-(m_k \mu_k - \lambda_k) t} - e^{-\mu_k t}\big),
\end{equation}
where, denoting by $E_{2, m_k}(r)$ the Erlang-C formula for the probability that a job has to wait in the M/M/$m_k$ queue with offered load $r$,
\begin{equation}
\label{eq:gklk}
G_k(\lambda_k)  \triangleq \frac{E_{2, m_k}\big(r_k(\lambda_k)\big)}{1-m_k + r_k(\lambda_k)}.
\end{equation}

Yet, note that (\ref{eq:wktlambdak}) is only valid when $r_k(\lambda_k) \neq m_k-1$, i.e., when $\lambda_k \neq  (m_k-1)\mu_k$. In the case $r_k(\lambda_k) = m_k-1$, we have
\begin{equation}
\label{eq:wktlambdak2}
S_k(t; \lambda_k)  = 
1 - \big(1 + \mu_k E_{2, m_k}(m_k-1)  t\big) e^{-\mu_k t}.
\end{equation}

Hence, the deadline-miss probability $\bar{P}_k^{\textup{miss}}(\lambda_k) \triangleq \Prob\{\widetilde{T}_k(\lambda_k) > \tau\}$ for a randomly arriving job that is admitted and assigned to queue $k$ is given by 
\begin{equation}
\label{eq:sbark}
\begin{split}
\bar{P}_k^{\textup{miss}}(\lambda_k) & = \int_0^\infty \big(1 - S_k(t; \lambda_k)\big) \, dF(t) \\
& = 
\begin{cases}
\phi(\mu_k) + G_k(\lambda_k)\big(\phi(m_k \mu_k - \lambda_k) - \phi(\mu_k)\big), & \lambda_k \neq  (m_k -1) \mu_k   \\
\phi(\mu_k) - \mu_k E_{2, m_k}(m_k-1)  \phi'(\mu_k), & \lambda_k = (m_k -1) \mu_k,
\end{cases}
\end{split}
\end{equation}
where $\phi(s) \triangleq \Ex\big[e^{-s \tau}\big] = \int_0^\infty e^{-s t} \, dF(t)$ is the Laplace transform of the  deadline $\tau$. We will consider the cases $\tau \equiv t$, with $\phi(s) = e^{-s t}$, $\tau \sim \mathrm{Unif}[t_{1}, t_{2}]$, with $\phi(s) = (e^{-t_{1} s} - e^{-t_{2} s})/(s (t_{2}-t_{1}))$, and $\tau \sim \textup{Exp}(\theta)$, with $\phi(s) = \theta/(s+\theta)$.



Thus, we can formulate the problem of finding the optimal BS as the following NLP, where $f_k(\lambda_k) \triangleq \lambda_k \bar{P}_k^{\textup{miss}}(\lambda_k)$ is the \emph{deadline-miss rate} for queue $k$:
\begin{equation}
\label{eq:nlp}
 \min \Bigg\{R \lambda_0 + \sum_{k \in \mathbb{K}} f_k(\lambda_k)\colon 
\sum_{k \in \mathbb{K}} \lambda_k = \lambda, \lambda_0
\geqslant 0,
0 \leqslant \lambda_k \leqslant m_k \mu_k, k \in \mathbb{K}\Bigg\}.
\end{equation}

Notice that, in problem (\ref{eq:nlp}), we allow a BS to saturate the queue at queue $k$ by allocating to it an arrival rate equal to its service capacity ($\lambda_k = m_k \mu_k$). In such a case ---which as we will see turns out to be relevant when $\lambda$ and $R$ are both large enough, and also in the pure-routing version of the problem, where  $\lambda_0 \equiv 0$, when $\lambda$ is large enough--- the probability that a random job assigned to a saturated queue $k$ misses its deadline is $\bar{P}_k^{\textup{miss}}(m_k \mu_k) = 1$, and hence we define $f_k(m_k \mu_k) \triangleq m_k \mu_k$.

\cite{kallmescass95} presents a Lagrangian algorithm for solving the version of (\ref{eq:nlp}) with strict inequalities $\lambda_k < m_k \mu_k$ in the  special case  $R = 1$ with  single-server queues $m_k = 1$, showing that taking $R = 1$ ensures that the optimal solution to (\ref{eq:nlp}) has $\lambda_k^* < m_k \mu_k$ for each $k$.
As for the multiserver case
$m_k \geqslant 2$, a Lagrangian algorithm is outlined  in \cite{heetal06}, albeit only for the pure-routing case where $\lambda_0  \equiv 0$, and without discussing the possibility that some queue(s) may be saturated in an optimal solution.

We next discuss how to extend the results in such works to obtain the optimal BS $\boldsymbol{\lambda}^* = (\lambda_k^*)_{k=0}^n$ solving (\ref{eq:nlp}) for the present model, assuming
that each function $\bar{P}_k^{\textup{miss}}(\lambda_k)$ satisfies the following condition:

\begin{itemize}
\item[\textup{(C)}] $\bar{P}_k^{\textup{miss}}(\lambda_k)$ is continuous, increasing and strictly convex on $[0, m_k \mu_k]$, having a continuous and increasing derivative  $(d/d\lambda_k)\bar{P}_k^{\textup{miss}}(\lambda_k)$ on $(0, m_k \mu_k),$ with finite one-sided derivatives 
$(d/d\lambda_k)\bar{P}_k^{\textup{miss}}(0^+)$ and 
$(d/d\lambda_k)\bar{P}_k^{\textup{miss}}\big((m_k \mu_k)^-\big)$.
$\alpha_k \triangleq f_k'(0^+)$ and $\beta_k \triangleq f_k'\big((m_k \mu_k)^-\big)$.
\end{itemize}

Under  (C), which is readily verified in the cases of constant, uniform and exponential deadlines considered here (e.g., in the case $\tau \equiv t$ for a single-server queue, $\bar{P}_k^{\textup{miss}}(\lambda_k) = e^{-(\mu_k - \lambda_k)t}$), 
the function $f_k(\lambda_k)$ inherits the same properties as $\bar{P}_k^{\textup{miss}}(\lambda_k)$, since $f_k'(\lambda_k) = \bar{P}_k^{\textup{miss}}(\lambda_k) + \lambda_k (d/d\lambda_k) \bar{P}_k^{\textup{miss}}(\lambda_k)$ on $(0, m_k \mu_k)$, having one-sided derivatives $\alpha_k \triangleq f_k'(0^+)$ and $\beta_k \triangleq f_k'\big((m_k \mu_k)^-\big)$ that satisfy 
\begin{equation}
\label{eq:osder}
0 < \alpha_k = \bar{P}_k^{\textup{miss}}(0)  = \phi(\mu_k) < 1 < 1 + m_k \mu_k (d/d \lambda_k)\bar{P}_k^{\textup{miss}}\big((m_k \mu_k)^-\big) = \beta_k.
\end{equation}
Hence, under (C), problem (\ref{eq:nlp}) has a unique optimal solution  $\boldsymbol{\lambda}^*$. 

Assume,  by reordering if necessary, that $\alpha_1 \leqslant \alpha_2 \leqslant \cdots \leqslant \alpha_n$. Note that, since $\alpha_k = \phi(\mu_k)$ and $\phi(s)$ is nonincreasing, such an ordering is equivalent to $\mu_1 \geqslant \mu_2 \geqslant \cdots \geqslant \mu_n$.
Further, let $x_1, \ldots, x_n$ be an ordering of the $n$ queues for which $\beta_{x_1} \leqslant \beta_{x_2} \leqslant \cdots \leqslant \beta_{x_n}$. 

It follows from the above  (cf.\ (\ref{eq:osder})) that, for any queue $k$ and  $\alpha \in [\alpha_k, \beta_k]$, the equation $f_k'(\lambda_k) = \alpha$ has a unique root $\Lambda_k^*(\alpha)$ in $[0 , m_k \mu_k]$, where $\Lambda_k^*(\alpha)$ is  continuous and increasing in $\alpha$,  and satisfies that $\Lambda_k^*(\alpha_k) = 0$ and $\Lambda_k^*(\beta_k) = m_k \mu_k$.
For convenience of notation, we also define 
\begin{equation}
\label{eq:lambdakstext}
\Lambda_k^*(\alpha) \triangleq  0, \quad\alpha < \alpha_k, \quad \text{and} \quad \Lambda_k^*(\alpha) \triangleq  m_k \mu_k, \quad \alpha > \beta_k.
\end{equation}

Now, by standard results in  convex optimization, the optimal BS solving  (\ref{eq:nlp}) is characterized by the following
 \emph{Karush--Kuhn--Tucker} (KKT) \emph{first-order optimality conditions}: a feasible BS $\boldsymbol{\lambda}^*$ for  (\ref{eq:nlp}) is optimal iff there exists a Lagrange multiplier $\alpha^*$ for the equality
constraint such that
\begin{equation}
\label{eq:kktc}
\begin{split}
R & \geqslant \alpha^*, \quad \text{ with ``=" if } \lambda_0^* > 0 \\
f_k'(\lambda_k^*)  & = \alpha^*, \quad \text{ for any queue $k$ with }  0 < \lambda_k^* < m_k \mu_k \\
\alpha_k  & \geqslant \alpha^*, \quad \text{ for any queue $k$ with } \lambda_k^* = 0 \\
 \beta_k  & \leqslant \alpha^*, \quad \text{ for any queue $k$ with }  \lambda_k^* = m_k \mu_k.
\end{split}
\end{equation} 

Note that it follows from (\ref{eq:kktc}) that, since $\beta_k > 1$, 
in the case $R = 1$ the optimal BS will have $\lambda_k^* < m_k \mu_k$ for every queue $k$.  

To  compute the  optimal BS $\boldsymbol{\lambda}^*$ for (\ref{eq:nlp}), we next present a ranking algorithm that
extend to the current setting the algorithm given in \citet[p.\ 321]{kallmescass95} for the special case $R = 1$ and $m_k \equiv 1$. 
The following algorithm determines an $\alpha^*$  such that the optimal BS satisfying (\ref{eq:kktc}) is given by $\lambda_k^* \triangleq \Lambda_k^*(\alpha^*)$ for $k \in \mathbb{K}$ ---taking into account (\ref{eq:lambdakstext})--- and 
$\lambda_0^* \triangleq \lambda - \sum_{k \in \mathbb{K}} \lambda_k^*$.
Three cases need to be distinguished (note that case 3 does not arise when $R = 1$):

\begin{description}
\item[Case 1 (light load):] for some $1 \leqslant l < n$, $\sum_{k=1}^l \Lambda_k^*(\alpha_l) < \lambda \leqslant \sum_{k=1}^l \Lambda_k^*(\alpha_{l+1})$. Letting $\widetilde{\alpha}$ be the unique root  in $(\alpha_l, \alpha_{l+1}]$  of $\sum_{k=1}^l \Lambda_k^*(\alpha) = \lambda$, take $\alpha^* \triangleq \min(\widetilde{\alpha}, R)$. Note that $\lambda_k^* >0$ for any $1 \leqslant k \leqslant l$ with $\alpha_k < R$, $\lambda_k^* = 0$ for any $1 \leqslant k \leqslant l$ with $\alpha_k \geqslant R$ and for $l < k \leqslant n$, and $\lambda_0^* = 0$ iff $\widetilde{\alpha} \leqslant R$.
\item[Case 2 (medium load):] $\sum_{k \in \mathbb{K}} \Lambda_k^*(\alpha_n) < \lambda < \sum_{k \in \mathbb{K}} \Lambda_k^*(\beta_{x_1})$. Letting $\widetilde{\alpha}$ be the unique root  in $(\alpha_n, \beta_{x_1})$ of $\sum_{k \in \mathbb{K}} \Lambda_k^*(\alpha) = \lambda$,  
take $\alpha^* \triangleq \min(\widetilde{\alpha}, R)$. Note that $0 < \lambda_k^* < m_k \mu_k$ for any $k \in \mathbb{K}$ with $\alpha_k < R$, $\lambda_k^* = 0$ for any $k \in \mathbb{K}$ with $\alpha_k \geqslant R$,  and $\lambda_0^* = 0$ iff $\widetilde{\alpha} \leqslant R$.
\item[Case 3 (heavy load):] for some $1 \leqslant l < n$, $\sum_{k \in \mathbb{K}} \Lambda_k^*(\beta_{x_l}) \leqslant \lambda < \sum_{k \in \mathbb{K}} \Lambda_k^*(\beta_{x_{l+1}})$. Letting $\widetilde{\alpha}$ be the unique root in $[\beta_{x_l}, \beta_{x_{l+1}})$ of $\sum_{k \in \mathbb{K}} \Lambda_k^*(\alpha) = \lambda$, take $\alpha^* \triangleq \min(\widetilde{\alpha}, R)$.
Note that $\lambda_{k}^* = m_k \mu_k$ for any $k \in \{x_1, \ldots, x_l\}$ with $\beta_k \leqslant R$,  $\lambda_{k}^* < m_k \mu_k$ for any $k \in \{x_1, \ldots, x_l\}$ with $\beta_k > R$ and for  $k \in \{x_{l+1}, \ldots, x_n\}$, $\lambda_{k}^* > 0$ for any $k \in \mathbb{K}$ with $\alpha_k < R$, and 
$\lambda_0^* = 0$ iff $\widetilde{\alpha} \leqslant R$.
\end{description}

Note that, in the light-load case, some queue(s) $k$ will not be used at all under the optimal BS. In contrast, in the heavy-load case some queue(s) $k$ may be saturated ($\lambda_k^* =  m_k \mu_k$). Note also that, since we assumed $\rho \triangleq \lambda/\sum_k m_k \mu_k < 1$, no higher arrival rates $\lambda$ need to be considered than those in Case 3.
Another observation is that, if the rejection cost is $R \geqslant \beta_{x_n}$ then no job is rejected under the optimal BS. At the opposite extreme, if $R \leqslant \alpha_1$ then all jobs are rejected.

\subsection{Second stage of the PI method: PI step}
\label{s:sspiim}
The second stage of the PI method is to carry out
one PI step of the policy iteration algorithm for MDPs, starting from the optimal BS policy, which yields a better policy.
Recall from Section \ref{s:toari} that an optimal policy, given by a solution to the DP optimality equation (\ref{eq:dpeqnsmodel}) for problem (\ref{eq:pa2}), takes the following form: upon arrival of a job finding the system  in state 
$\mathbf{x}$,  (i) reject it if $R \leqslant \min_k \Delta_k \widetilde{b}^*(\mathbf{x} + \mathbf{e}_k)$; otherwise, 
(ii) send it to a queue with lowest
   $\Delta_k \widetilde{b}^*(\mathbf{x} + \mathbf{e}_k)$.

The PI method replaces the relative cost increments $\Delta_k \widetilde{b}^*(\mathbf{x})$ 
by those under the optimal BS $\boldsymbol{\lambda}^*$, which we denote by $\Delta_k \widetilde{b}(\mathbf{x}; \boldsymbol{\lambda}^*)$. Now, 
since under a BS policy the $n$ queues evolve independently, separately accruing costs  in the modified problem (\ref{eq:pa2}) at rate $\lambda R  a_k(t) + h_k\big(X_k(t)\big)$ for queue $k$, we have $\widetilde{b}(\mathbf{x}; \boldsymbol{\lambda}^*) \triangleq \sum_k \widetilde{b}_k(x_k; \lambda_k^*)$, where 
$\widetilde{b}_k^*(x_k; \lambda_k^*)$ is the relative cost function for  queue $k$, considered in isolation.
Hence, $\Delta_k \widetilde{b}(\mathbf{x}; \boldsymbol{\lambda}^*) = \Delta b_k(x_k; \lambda_k^*)
\triangleq \widetilde{b}_k(x_k; \lambda_k^*) - \widetilde{b}_k(x_k-1; \lambda_k^*)$. The resulting PI policy is:
upon arrival of a job finding the system  in state 
$\mathbf{x}$, (i) reject it if $R \leqslant \min_k \Delta \widetilde{b}_k(x_k + 1; \lambda_k^*)$; otherwise,
(ii) route it to a queue $k$ with lowest
   $\Delta \widetilde{b}_k(x_k + 1; \lambda_k^*)$.

Hence, the PI policy is an index policy, having routing index
\begin{equation}
\label{eq:nukpii}
\nu_k^{\textup{PI}}(x_k) \triangleq \Delta \widetilde{b}_k(x_k+1; \lambda_k^*),
\end{equation}
for queue $k$, and global admission index $\gamma^{\textup{PI}}(\mathbf{x}) \triangleq \min_k \nu_k^{\textup{PI}}(x_k)$.

To evaluate the PI index $\nu_k^{\textup{PI}}(x_k)$ for queue $k$, 
we must solve the corresponding \emph{Poisson equations},  which are readily reformulated in terms of the index as
\begin{equation}
\label{eq:pepiim}
\begin{split}
f_k(\lambda_k^*) - \lambda_k^* \nu_k^{\textup{PI}}(0)  & = h_k(0) \\
f_k(\lambda_k^*) - \lambda_k^* \nu_k^{\textup{PI}}(x_k) + \bar{\mu}_k(x_k)
\nu_k^{\textup{PI}}(x_k-1) & = h_k(x_k), \quad  x_k \geqslant 1,
\end{split}
\end{equation}
where $\bar{\mu}_k(x_k) \triangleq \mu_k \min(x_k, m_k)$, and  $h_k(x_k)$ is the holding cost rate defined in (\ref{eq:ridef}). 
Thus, (\ref{eq:pepiim}) gives a first-order linear recursion for computing the  index values
$\nu_k^{\textup{PI}}(0), \ldots, \nu_k^{\textup{PI}}(x_k)$ in linear, $O(x_k)$ time, given $\lambda_k^*$ and $f_k(\lambda_k^*)$. 

We write below $\rho_k^* \triangleq \rho_k(\lambda_k^*)$, $\bar{P}_k^{\textup{miss}, *} \triangleq  \bar{P}_k^{\textup{miss}}(\lambda_k^*)$ and 
$f_k^* \triangleq  f_k(\lambda_k^*)$. Note that $f_k^* = \lambda_k^* \bar{P}_k^{\textup{miss}, *}$.

\begin{figure}[!ht]
\centering
\includegraphics[height=2in]{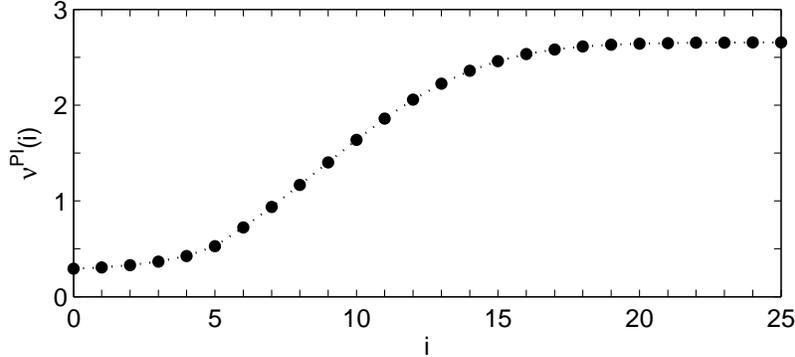}
\caption{The PI index for a queue with $m = 6$, $\mu = 0.5$, $\tau \sim \mathrm{Unif}[2.5, 3.5]$, and $\lambda^* = 2.1$.}
\label{fig:plotPIindex}
\end{figure}

\begin{remark}
\label{re:frro} {\rm 
The following results are readily obtained.
\begin{enumerate}
\item In the case of a queue $k$ that is not used under the optimal BS, i.e., having $\lambda_k^* = 0$, (\ref{eq:pepiim}) and (\ref{eq:ridef}) yield that $\nu_k^{\textup{PI}}(x_k) = P_k^{\textup{miss}}(x_k)$, and hence the PI index reduces to the IO index $\nu_k^{\textup{IO}}(x_k)$.
\item In the case $0 < \lambda_k^* < m_k \mu_k$, i.e., $0 < \rho_k^* < 1$, $\nu_k^{\textup{PI}}(0) = \bar{P}_k^{\textup{miss}, *}$ and
\begin{equation}
\label{eq:nukssq}
\nu_k^{\textup{PI}}(x_k) = \frac{1}{(\rho_k^*)^{x_k}}\bigg[\nu_k^{\textup{PI}}(m_k-1) + \sum_{j_k=m_k-1}^{x_k-1} (\rho_k^*)^{j_k} \big(\rho_k^* \bar{P}_k^{\textup{miss}, *} - P_k^{\text{miss}}(j_k)\big)\bigg], \, x_k \geqslant m_k.
\end{equation}
\item In the case of a saturated queue ($\rho_k^* = 1$), $\bar{P}_k^{\textup{miss}, *} = 1$, $\nu_k^{\textup{PI}}(0) = 1$, and
\begin{equation}
\label{eq:nukssq2}
\nu_k^{\textup{PI}}(x_k) = \nu_k^{\textup{PI}}(m_k-1) + \sum_{j_k=m_k-1}^{x_k-1} \big(1 - P_k^{\text{miss}}(j_k)\big), \, x_k \geqslant m_k.
\end{equation}
\item If $\nu_k^{\textup{PI}}(x_k)$ converges as $x_k\to \infty$ to a finite limit $\nu_k^{\textup{PI}}(\infty) $, then 
\begin{equation}
\label{eq:nuklim}
\nu_k^{\textup{PI}}(\infty) = 
\begin{cases}
\displaystyle \frac{m_k \mu_k -f_k^*}{ m_k \mu_k - \lambda_k^*} = \frac{1 - \rho_k^* \bar{P}_k^{\textup{miss}, *}}{1-\rho_k^*} > 1, & \text{ if } 0 < \rho_k^* < 1 \\
\displaystyle \nu_k^{\textup{PI}}(m_k-1) + \sum_{j_k=m_k-1}^{\infty} \big(1 - P_k^{\text{miss}}(j_k)\big), & \text{ if } \rho_k^* = 1.
\end{cases}
\end{equation}
Based on extensive numerical experience, we conjecture that the PI index $\nu_k^{\textup{PI}}(x_k)$ for the present model is monotone increasing and converges to a finite limit. See, e.g., Figure \ref{fig:plotPIindex}.
Further, L'H\^opital's rule yields that $\lim_{\lambda_k^* \nearrow m_k \mu_k} \big(m_k \mu_k -f_k(\lambda_k^*)\big)/\big(m_k \mu_k - \lambda_k^*\big) = \beta_k \triangleq f_k'\big((m_k \mu_k)^-\big)$, suggesting that the limiting index in the case $\rho_k^* = 1$ should be $\nu_k^{\textup{PI}}(\infty) =\beta_k$, which the author has also verified numerically.
Yet, such experiments also reveal that the computed index $\widehat{\nu}_k^{\textup{PI}}(x_k)$  may occasionally diverge in practice, due to small inaccuracies in the computed value of $f_k^*$. We propose to overcome such numerical instability  by setting the computed index value equal to the appropriate limit in (\ref{eq:nuklim}) as soon as $\widehat{\nu}_k^{\textup{PI}}(x_k)$ shows any signs of diverging, e.g., by exceeding the limit $\nu_k^{\textup{PI}}(\infty)$ in (\ref{eq:nuklim}).
\end{enumerate}}
\end{remark}

\section{RB index}
\label{s:mod2mp}
The RB routing index for a queue is obtained by adapting to the 
present model the approach and results introduced in
\cite{nmmp02} for the design of policies for control of admission and routing to parallel queues via  Whittle's RB indexation (see \cite{whit88}). See also \cite{nmnetcoop07}.
To define Whittle's RB index, we consider the single-queue case ($n = 1$) of  modified problem (\ref{eq:pa2}), concerning the optimal dynamic control of admission to an M/M/$m$ queue with arrival rate $\lambda$ and service rate $\mu$ at each server. Note that we drop the queue label $k$. Jobs come with a response-time deadline  with CDF $F(t)$. At the arrival instant of a job, the controller can choose to reject it, incurring a rejection penalty of $R$, or accept it. In the latter case, a unit penalty is incurred if the job's deadline is missed. As discussed in Section \ref{s:rjosor}, we will account for deadline-miss costs via continuously accruing action-independent holding costs, which are incurred at rate $h(x)$ (given by (\ref{eq:ridef})) per unit time when there are $x$ jobs in the queue.

Consider the corresponding \emph{optimal admission control problem}, which is to find an admission control policy $\pi^*$, drawn from the class $\Pi$ of policies that base decisions on the history of queue states and actions, that minimizes the average rate of rejection and deadline-miss costs  per unit time. 
Adapting the notation in (\ref{eq:pa2}) to the present setting,  we  formulate such a problem as
\begin{equation}
\label{eq:pasqh}
\min_{\pi \in \Pi} \limsup_{T \to \infty} \frac{1}{T}
\Ex_{i^0}^{\pi}\Bigg[\int_0^T  \Big\{R \lambda  a(t) + h\big(X(t)\big)\Big\} \, dt\Bigg].
\end{equation}

By standard results, problem (\ref{eq:pasqh}) has an optimal policy that is stationary deterministic  and independent of the initial state.
Now, instead of considering a fixed rejection charge $R$, we will view it as  a \emph{parameter} that is allowed to take any real value $R \in \mathbb{R}$. 
Rather than solving each individual problem (\ref{eq:pasqh}) separately for different values of $R$, we will characterize optimal policies for the entire parametric family of such \emph{$R$-charge problems}, as $R \in \mathbb{R}$, in terms of a numeric index  attached to queue states.

Let us say that problem (\ref{eq:pasqh}) is \emph{indexable} if there exists a numeric index $\nu^{\textup{RB}}(x)$ attached to queue states $x$, which does not depend on $R$,   such that, for any rejection charge $R \in \mathbb{R}$ and  state $x$,  it is optimal to reject an arrival when the queue holds $x$ jobs iff $\nu^{\textup{RB}}(x) \geqslant R$.

Note that, in particular, if problem (\ref{eq:pasqh}) is indexable and the index $\nu^{\textup{RB}}(x)$ is nondecreasing in the state $x$, then it follows that problem (\ref{eq:pasqh}) is solved optimally by a \emph{threshold policy} (reject an arrival iff the queue length is larger than a certain threshold). 
While the conventional approach to establish the optimality of threshold policies for admission control problems such as (\ref{eq:pasqh}) relies on assuming convexity of the holding cost rate function $h(x)$ (see \cite{stid85}), we will be able to establish the stronger result of \emph{indexability} with a nondecreasing index $\nu^{\textup{RB}}(x)$, even though the holding cost rate $h(x)$ defined in (\ref{eq:ridef}) is \emph{not convex} (note that $h(x) = m \mu P^{\textup{miss}}(x-1)$ for $x \geqslant m$, and hence $h(x)$ grows to $m \mu$ as $x \to \infty$).

To establish that problem (\ref{eq:pasqh}) is indexable with a monotone nondecreasing RB index $\nu^{\textup{RB}}(x)$, and to evaluate the index, we will deploy the general sufficient conditions for indexability and the \emph{adaptive-greedy index algorithm} introduced in \cite{nmaap01,nmmp02}, which do not require such convexity assumptions.
 Such conditions are formulated in  \citet{nmmp02} in terms of certain quantities $w^{S_y}(x)$, for $x \geqslant 0, y \geq 1$,  which represent \emph{marginal rejection measures}, and 
 indices $\nu^{\textup{RB}}(x)$, which are produced by the algorithm. 
Here, we will write $w(x) \triangleq w^{S_{x+2}}(x)$.

In the present model, 
such an index algorithm reduces to a coupled set of first-order linear recursions (see \citet[pp.\ 396--397]{nmmp02}). First, compute
\begin{equation}
\label{eq:wrcomp}
\begin{split}
z(1) & = 1; \quad z(x) =  1-\frac{\lambda  \bar{\mu}({x-1})}{\big[\lambda +
\bar{\mu}({x-1})\big]  \big[\lambda + \bar{\mu}(x)\big] z({x-1})}, \, x \geqslant 2, \\
w(0) & = 
\frac{\lambda \mu}{\lambda + \mu}; \quad
w({x-1})  = \lambda
\frac{\displaystyle \Delta \bar{\mu}(x) + 
  \frac{w({x-2})}{\bar{\rho}({x-2})}}{\big[\lambda + \bar{\mu}(x)\big] z(x)}, \, x \geqslant 2,
\end{split}
\end{equation}
where $\bar{\rho}(x) \triangleq
\lambda/\bar{\mu}({x+1})$, $\bar{\mu}(x) \triangleq \min(x, m) \mu$, and $\Delta \bar{\mu}(x) \triangleq \bar{\mu}(x) - \bar{\mu}(x-1)$. Then, recursively compute the  index $\nu^{\textup{RB}}(x)$ by
\begin{equation}
\label{eq:windexqs}
\begin{split} 
\nu^{\textup{RB}}(0) & = \frac{h(1)}{\mu}; \,
\nu^{\textup{RB}}(x)  = 
\nu^{\textup{RB}}(x-1) + \frac{\displaystyle{\Delta h(x+1) - 
  \nu^{\textup{RB}}(x-1)} \Delta \bar{\mu}(x+1)}{\displaystyle{\Delta \bar{\mu}(x+1) + 
  \frac{w(x-1)}{\bar{\rho}(x-1)}}}.
\end{split}
\end{equation}

We can simplify (\ref{eq:wrcomp})--(\ref{eq:windexqs}) by noting that $\Delta \bar{\mu}(x) = \mu$ for $1 \leqslant x \leqslant m$, and
$\Delta \bar{\mu}(x) = 0$ for $x \geqslant m+1$. Further, from (\ref{eq:ridef}) and $P^{\textup{miss}}(x) = P^{\textup{miss}}(0)$ for $1 \leqslant x < m$ (see Section \ref{s:model1io}), we have $\Delta h(x) \triangleq  h(x) - h(x-1) =
 \mu P^{\textup{miss}}(0)$ for $i = 1, \ldots, m$, whereas $\Delta h(x) = m \mu \Delta P^{\textup{miss}}(x-1)$ for $x \geqslant m+1$. Hence, 
\begin{equation}
\label{eq:nurbi1m}
\nu^{\textup{RB}}(x) = 
\begin{cases}
P^{\textup{miss}}(0),& 0 \leqslant  x \leqslant m-1 \\
\displaystyle \nu^{\textup{RB}}(x-1) + \lambda \frac{\Delta P^{\textup{miss}}(x)}{
  w(x-1)}, &  x \geqslant m.
\end{cases}
\end{equation}

Such a set of first-order linear recursions, along with those in Section \ref{s:model1io} for $P^{\textup{miss}}(x)$ (which equals $\nu^{\textup{IO}}(x)$), efficiently computes the index values $\nu^{\textup{RB}}(0)$, \ldots,
$\nu^{\textup{RB}}(x)$ in $O(x)$ time.

In the single-server case $m = 1$ such recursions are readily solved:
\[
z(x) = 
\begin{cases}
\displaystyle \frac{\rho ^{x+1}-1}{(\rho +1) \left(\rho ^x-1\right)} & \text{ if } \rho \neq 1
\\ 
\displaystyle \frac{x+1}{2 x} & \text{ if } \rho = 1
\end{cases}
\quad \text{and} \quad
w(x) = \begin{cases}
\displaystyle \lambda  \frac{\rho -1}{\rho ^{x+2}-1}
& \text{ if } \rho \neq 1
\\
\displaystyle 
   \frac{\lambda}{x+2} & \text{ if } \rho = 1.
\end{cases}
\]
Further, assuming deterministic deadlines $\tau \equiv t$, we have
$\Delta P^{\textup{miss}}(x) = e^{-\mu t} (\mu t)^x/x!$. 
Hence, in such a case the RB index is given by $\nu^{\textup{RB}}(0) = e^{-\mu t}$, and, for $x \geqslant 1$,
\begin{equation}
\label{eq:nurbimeq1}
\nu^{\textup{RB}}(x) = \begin{cases}
\displaystyle \frac{e^{-\mu t}}{\rho-1} \sum_{j=0}^x \frac{\rho (\lambda t)^j - (\mu t)^j}{j!} & \text{ if } \rho \neq 1 \\
\displaystyle 
  e^{-\mu t} \sum_{j=0}^x \frac{(j+1)(\mu t)^j}{j!} & \text{ if } \rho = 1.
\end{cases}
\end{equation}

 The sufficient conditions for indexability referred to above, as they apply to the present model, are:
 \begin{itemize}
 \item[(i)] The marginal rejection measures  $w^{S_y}(x)$ are positive; and
  \item[(ii)] the index $\nu^{\textup{RB}}(x)$ is monotone nondecreasing in $x$.
 \end{itemize}
 
 The following result is established in \citet[Cor.\ 2]{nmaap01} and \citet[Th.\ 6.3]{nmmp02}, at increasing levels of generality. 
 \begin{theorem}
 \label{the:sicpcl}
Under conditions \textup{(i)--(ii)} the model is indexable, with index  $\nu^{\textup{RB}}(x)$.
 \end{theorem}
 
We are now ready to establish that the present model is indexable. 
\begin{proposition}
\label{pro:pmindx}
Problem $(\ref{eq:pasqh})$ is indexable, and $\nu^{\textup{RB}}(x)$ is its RB index.
\end{proposition}
\begin{proof}
Regarding  condition (i), its satisfaction is ensured by \citet[Prop.\ 7.2]{nmmp02} when service rates $\bar{\mu}(x)$ are concave nondecreasing in $x$, as is the case here.
As for condition (ii), it follows from (\ref{eq:nurbi1m}), noting that $P^{\textup{miss}}(x)$ is nondecreasing and $w(x-1) > 0$. Invoking Theorem \ref{the:sicpcl} completes the proof.
\end{proof}

\begin{figure}[!ht]
\centering
\includegraphics[height=2in]{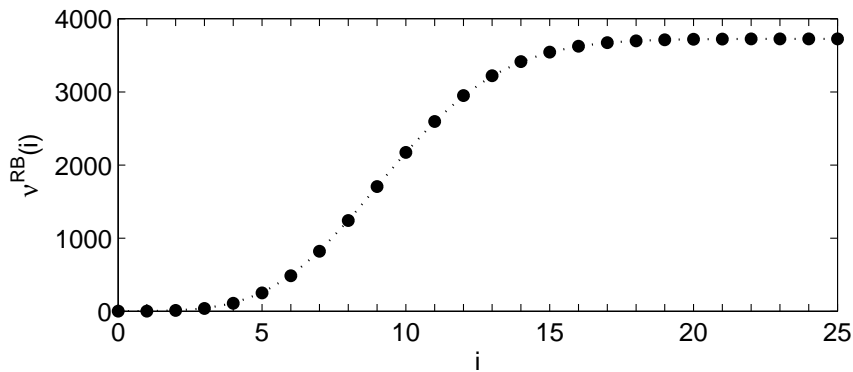}
\caption{The RB index for a queue with $m = 1$, $\mu = 1$, $\tau \equiv 2$, and $\lambda = 5$.}
\label{fig:plotRBindex}
\end{figure}

It must be emphasized that, in the system model of concern in this paper, with $n$ parallel queues, 
the above RB index is to be computed separately for each of the queues, using for queue $k$  in the above formulae the overall arrival rate $\lambda$ and the parameters $m_k$ and $\mu_k$. Hence, unless the parallel-queue system is lightly loaded, one will typically have $\rho_k \triangleq \lambda/(m_k \mu_k) > 1$ for some or even all queues $k$, even while the system load  $\lambda/\sum_k m_k \mu_k$ remains below unity. 
In such a typical case, 
the computation and use of the RB index raises the following issues:
\begin{enumerate}
\item It is readily verified that, for a queue with $\lambda > m \mu$, i.e., $\rho > 1$, $\lim_{x \to \infty}  w(x) = 0$. Since $\lim_{x \to \infty}  \Delta P^{\textup{miss}}(x)= 0$, computation of the RB index via (\ref{eq:nurbi1m}) can be numerically unstable.
\item Numerical experiments reveal that the RB index $\nu^{\textup{RB}}(x)$ grows extremely fast with $\lambda$, converging to a finit limit $\nu^{\textup{RB}}(\infty)$ as $x \to \infty$, which can be very large for a queue with utilization $\rho > 1$. 
The cause of this phenomenon is apparent in the single-server case $m = 1$ with a deterministic deadline $\tau \equiv t$, in which (\ref{eq:nurbimeq1}) yields that the index $\nu^{\textup{RB}}(x)$ converges as $x \to \infty$ to the limit
\begin{equation}
\label{eq:nurblimit}
\nu^{\textup{RB}}(\infty) = \begin{cases}
\displaystyle \frac{\rho e^{(\lambda - \mu) t} - 1}{\rho-1} & \text{ if } \rho \neq 1 \\
\displaystyle 
  1 + \mu t & \text{ if } \rho = 1.
\end{cases}
\end{equation}
Hence, $\nu^{\textup{RB}}(\infty)$ grows exponentially with $\lambda$. Thus, e.g., for a queue with $m = 1$, $\mu = 1$, $\tau \equiv 2$, and $\lambda = 5$, $\nu^{\textup{RB}}(x)$ converges as $x \to \infty$ to $\nu^{\textup{RB}}(\infty) = 3725.95$, as shown in Figure  \ref{fig:plotRBindex}.
As a consequence, when used to design an admission policy (reject an arrival if the rejection charge $R$ is lower than the current  index value at every queue), the  resulting admission policy may grossly overreject jobs relative to what is optimal.
\item The convergence of the index $\nu^{\textup{RB}}(x)$ to the finite limit $\nu^{\textup{RB}}(\infty)$ implies that, in the single-queue admission control problem (\ref{eq:pasqh}), it is optimal to admit \emph{all} arrivals iff the rejection cost per job $R$ is at or above $\nu^{\textup{RB}}(\infty)$.
\end{enumerate}

\section{Numerical experiments}
\label{s:ne}
This section reports on the results of a numerical study, whose aim is to benchmark the cost performance of the policies considered above (the static BS policy, and the dynamic IO, PI, and RB index policies), both against each other and against the optimal cost performance, in a variety of scenarios.

To  accurately evaluate the optimal cost performance, attention was restricted to instances with $n = 2$ queues (clusters). 
The optimal average cost per job was evaluated by solving with CPLEX the linear programming formulation of the average-cost DP equations (\ref{eq:dpeqnsmodel}), truncating the buffer size of each queue to 60 jobs (it was verified that increasing buffer sizes gave the same results for the instances considered).
As for the performance of each of the dynamic index policies, it was evaluated by solving with MATLAB the corresponding Poisson equations.

The study encompasses both the model of optimal joint control of admission and routing, and the pure-routing model where the admission control capability is disabled.

The experiments were designed to assess the effect, on the performance of each policy, of varying the system load, the rejection cost, the degree of heterogeneity in either cluster server speeds or server-pool sizes, and the magnitude and distribution of the response-time deadline.
Recall from Section \ref{s:intro} that the deadline-miss cost is normalized to the value $C = 1$. Further, in light of the alternative interpretation given in Section \ref{s:intro} of $R$ as the admission fee per job and $C$ as the deadline-miss refund, we only consider instances where $R \geqslant C$, i.e., where the refund per missed deadline is a fraction of the admission fee, as this appears more reasonable in practice than the case $R < C$.
Note that, in such cases, the IO index policy doest not reject any jobs, since its value is below $1$.

\subsection{Effect of system load}
\label{s:esl}
\begin{figure}[!htb]
\centering
\includegraphics[height=3.5in]{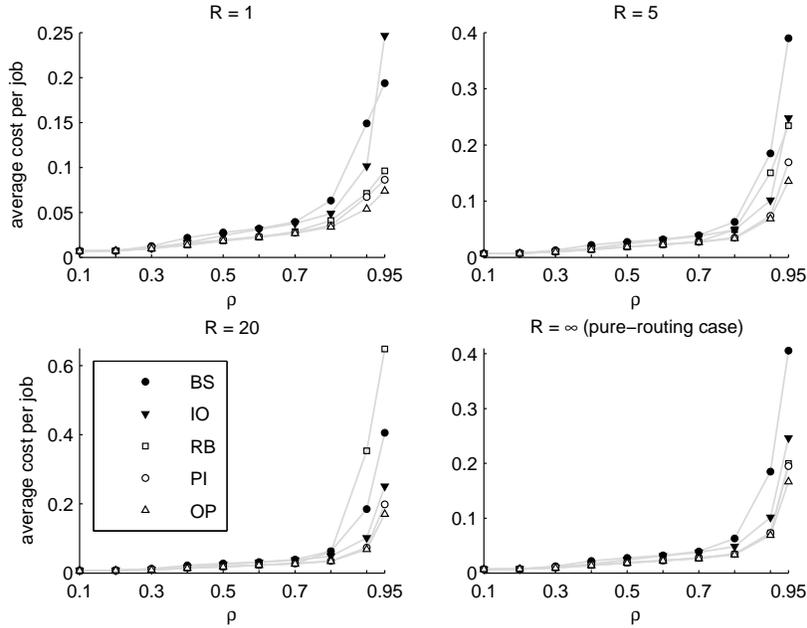}
\caption{Effect of system load with $\tau \equiv 1$.}
\label{fig:exp1a}
\end{figure}

We consider a base instance having $n = 2$ queues, with server-pool sizes $\mathbf{m} = (4, 8)$ and service rates $\boldsymbol{\mu} = (5, 3)$. The response-time deadline is deterministic, being $\tau \equiv 1$.
Figure  \ref{fig:exp1a} plots, for each of four values of the rejection cost $R \in \{1, 5, 20, \infty\}$, the average cost per job under each policy, and under an optimal policy, against the system load $\rho \triangleq \lambda/\sum_k m_k \mu_k$, which ranges over $\rho \in [0.1, 0.95]$ by changing the arrival rate $\lambda$.

The results are displayed in Figure \ref{fig:exp1a}, which shows that, while all the policies considered are nearly optimal under light loads, they exhibit significant differences in their cost performances under heavier loads.
For the model with admission control, the 
first three panes  show that 
the PI policy  consistently achieves a near-optimal performance.
The performance of the RB policy steeply deteriorates as the load gets heavier in cases with high rejection costs, becoming the worst policy for $\rho = 0.9$ and $R = 20$. Recall that, as noted above, the RB index can reach extremely high values under heavy loads, and hence yields admission policies that overreject.
As for the IO policy, its optimality loss worsens under as the load gets heavier in cases with low rejection costs (in contrast to the RB policy). This is to be expected, since in such cases it will be optimal to reject some jobs, whereas the IO policy does not reject any job.
The optimal BS policy steadily worsens as the load gets higher.

\begin{figure}[!htb]
\centering
\includegraphics[height=3.5in]{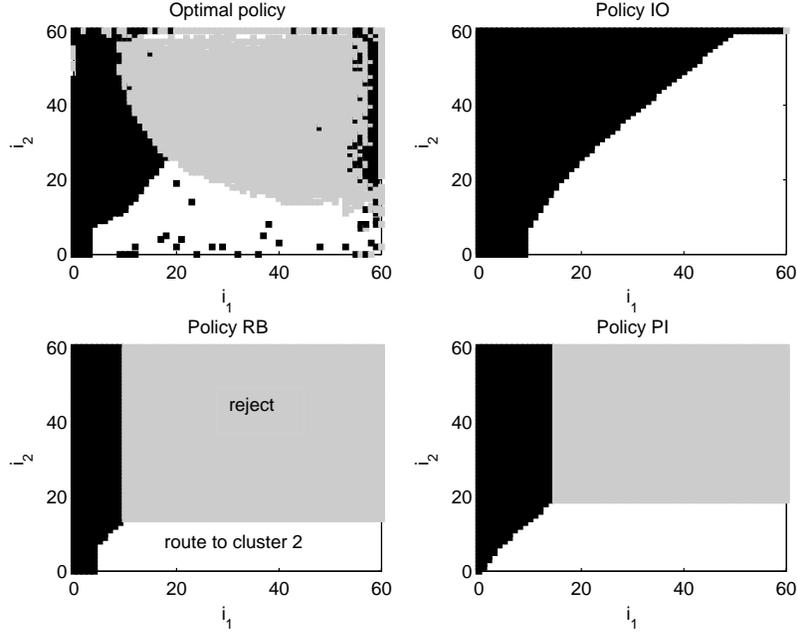}
\caption{Structure of computed policies for the case $\tau \equiv 1$, $\rho = 0.9$, $R = 4$.}
\label{fig:policyCompare}
\end{figure}

Figure \ref{fig:policyCompare} plots the structure of the index policies and the optimal policy that was computed in one of the instances. It is apparent that the RB policy rejects too many jobs. Although the shape of the rejection region differs under the PI policy and under an optimal policy, the former's performance is near optimal. Note that, in the plot of the optimal policy, the actions at some joint states $(x_1, x_2)$ differ from the overall pattern (e.g., there are some black dots inside the white region). This is due to numerical effects, as the linear programs from which such a solution results are rather large. In fact, we have found that in most instances the structure of the optimal policy is obscured due to such effects.

As for the pure-routing model, the lower right pane in Figure \ref{fig:exp1a} shows that both the PI and the RB policies are nearly optimal for the range of  loads considered. The IO policy deviates more from the optimal performance in heavy traffic. The performance of the BS policy   worsens substantially as traffic gets heavier.

\begin{figure}[!htb]
\centering
\includegraphics[height=3.8in]{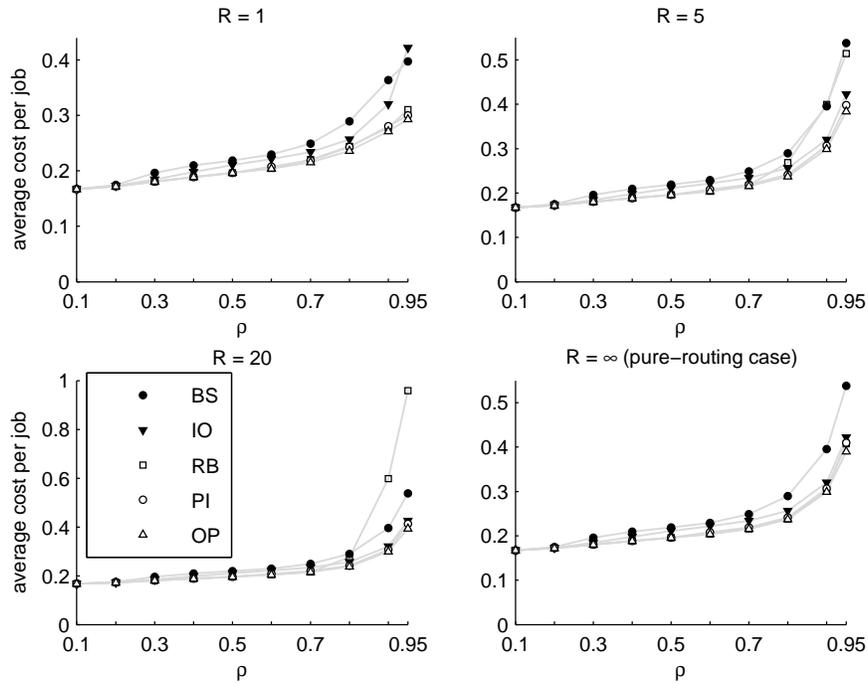}
\caption{Effect of system load with $\tau \sim \textup{Exp}(1)$.}
\label{fig:exp5a}
\end{figure}

\begin{figure}[!htb]
\centering
\includegraphics[height=3.8in]{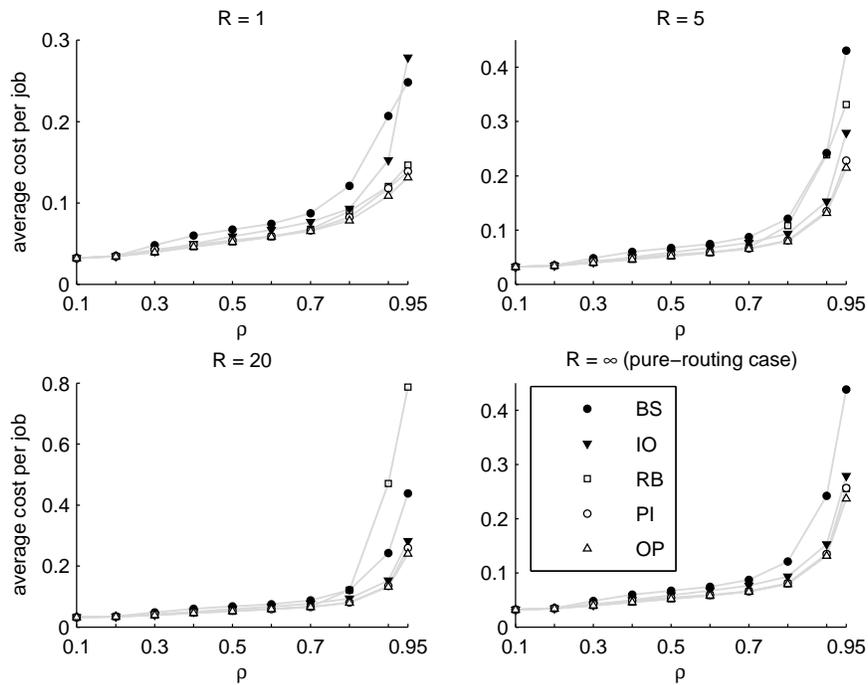}
\caption{Effect of system load with $\tau \sim \textup{Uniform}[0.3, 1.7]$.}
\label{fig:exp5b}
\end{figure}

Figures  \ref{fig:exp5a} and \ref{fig:exp5b} show the results of carrying out this experiment for two nondeterministic cases of  the  response-time deadline distribution: $\tau \sim \textup{Exp}(1)$ and  
$\tau \sim \textup{Uniform}[0.3, 1.7]$. Note that, in each  case, the mean response-time deadline is $1$. 
The results are qualitatively similar to those for a deterministic response-time deadline $\tau \equiv 1$.
Yet, notice that the level of average costs increase with the variability of the response-time deadline distribution, being highest in the exponential case, and lowest in the deterministic case.


\subsection{Effect of rejection cost}
\label{s:erc}
\begin{figure}[!ht]
\centering
\includegraphics[height=4in]{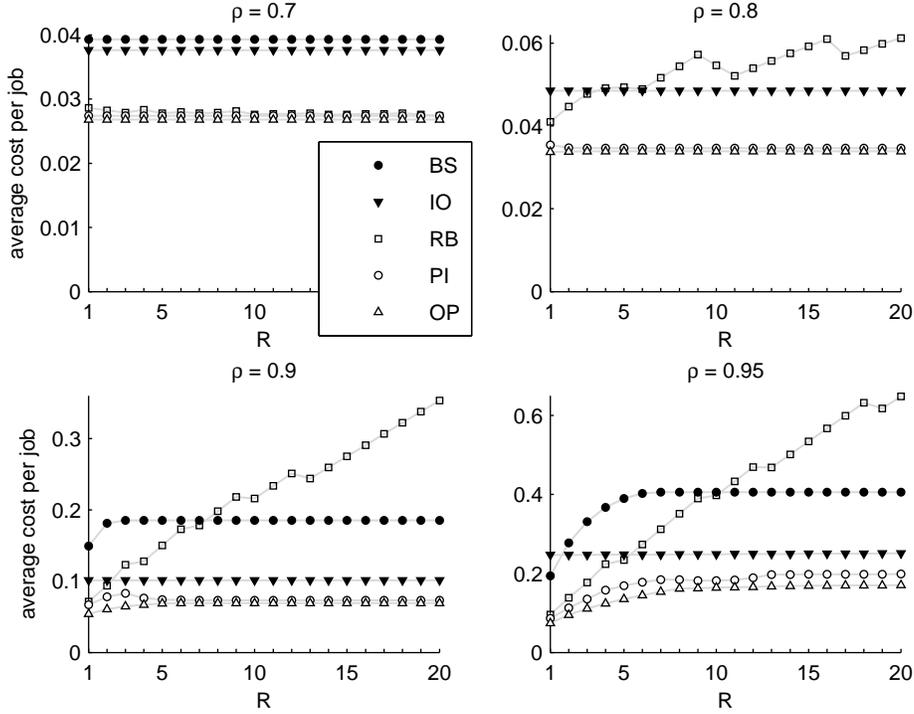}
\caption{Effect of rejection cost.}
\label{fig:exp1b}
\end{figure}

The second experiment aims to benchmark the performance of the policies considered as the rejection cost $R$ varies. 
We consider the same  base instance as in the first experiment. 
Figure  \ref{fig:exp1b} plots, for each of four values of the system load $\rho \in \{0.7, 0.8, 0.9, 0.95\}$, the average cost per job obtained under each policy, and under an optimal policy, against the rejection cost $R \in \{1, 2, \ldots, 20\}$. 

The results show that the performance of all policies considered, except RB, is not substantially affected by the rejection cost $R$. Such is also the case for the RB policy when the system load is low, but, as the load gets higher, its performance severely degrades as $R$ grows. This is due to the fact that, for high loads, the RB index takes very large values, and hence the resulting admission policy rejects too many jobs.

\subsection{Effect of heterogeneity in cluster server speeds} 
\label{s:ess}
\begin{figure}[!ht]
\centering
\includegraphics[height=4in]{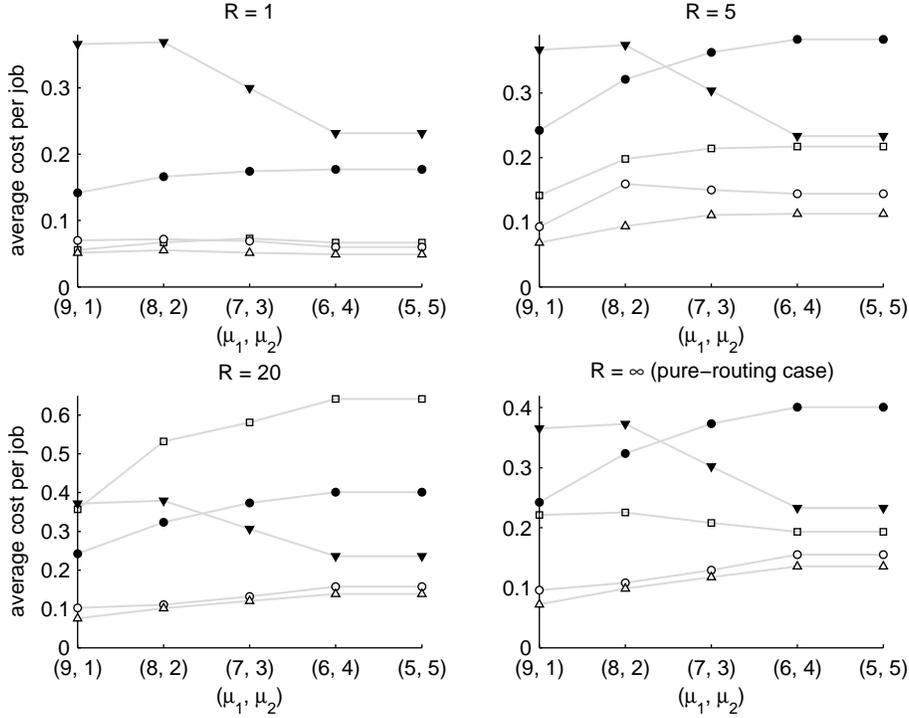}
\caption{Effect of heterogeneity in cluster server speeds.}
\label{fig:exp2}
\end{figure}

The third experiment aims to compare the policies as server speeds at the two queues range from a highly heterogeneous system to a homogeneous system. 
We consider a base instance with homogeneous server-pool sizes $\boldsymbol{m} = (4, 4)$, system load $\rho = 0.9$, and response-time deadline $\tau \equiv 1$.
Figure  \ref{fig:exp2} plots, for each of four values of the rejection cost $R \in \{1, 5, 20, \infty\}$, the average cost per job under each policy, and under an optimal policy, against five combinations of  service rates: $\boldsymbol{\mu} \in \{(9, 1), (8, 2), (7, 3), (6, 4), (5, 5)\}$. Note that the average speed per server the same in each combination.

For the model with admission control, the 
first three panes in Figure \ref{fig:exp2} show that
the cost performances of both the BS and the 
 RB policies deteriorate as cluster server speeds get more homogeneous, with BS (resp.\ RB) being the worse of the two for lower (resp.\ higher) rejection costs.
As for the IO policy, its optimality loss displays the opposite behavior, becoming closer to optimal as cluster server speeds get more homogeneous. Note that IO is the worst policy for $R = 1$.
The PI policy is, again, consistently near optimal.

As for the pure-routing model, the lower right pane in Figure \ref{fig:exp2} shows that the optimality loss of the BS policy deteriorates substantially as the queues become more homogeneous. Policies RB and IO display the opposite behavior, becoming closer to optimal as the clusters get more homogeneous. The best of the index policies if PI, followed by RB and then IO.

\subsection{Effect of heterogeneity in cluster server-pool sizes}
\label{s:esps}
\begin{figure}[!ht]
\centering
\includegraphics[height=4in]{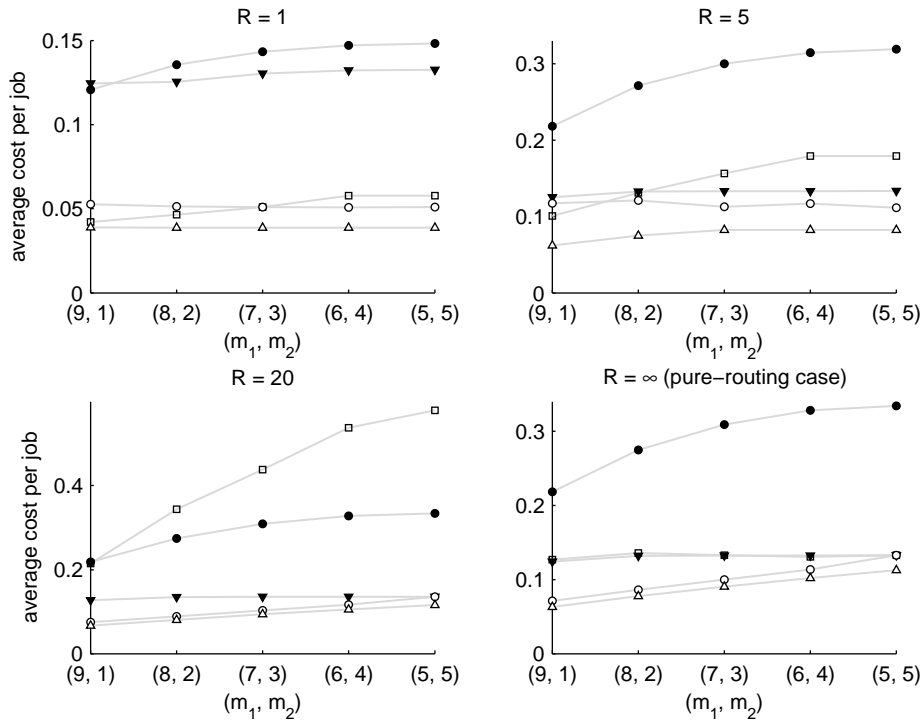}
\caption{Effect of heterogeneity in cluster  server-pool sizes.}
\label{fig:exp3}
\end{figure}

The aim of the fourth experiment is to benchmark the performance of the policies as server-pool sizes at the two queues range from a highly heterogeneous system to a homogeneous system. 
We consider a base instance with homogeneous server speeds $\boldsymbol{\mu} = (5, 5)$, system load $\rho = 0.9$, and response-time deadline $\tau \equiv 1$.
Figure  \ref{fig:exp3} plots, for each of four values of the rejection cost $R \in \{1, 5, 20, \infty\}$, the average cost per job under each policy, and under an optimal policy, against five combinations of  server-pool sizes: $\boldsymbol{m} \in \{(9, 1), (8, 2), (7, 3), (6, 4), (5, 5)\}$. Note that the total number of servers is the same in each combination.
The results of this experiment are qualitatively similar to those of the previous experiment in terms of ranking of the policies. Again, PI is the only policy that is consistently near optimal.

\subsection{Effect of magnitude of response-time}
\label{s:esrtg}
\begin{figure}[!ht]
\centering
\includegraphics[height=4in]{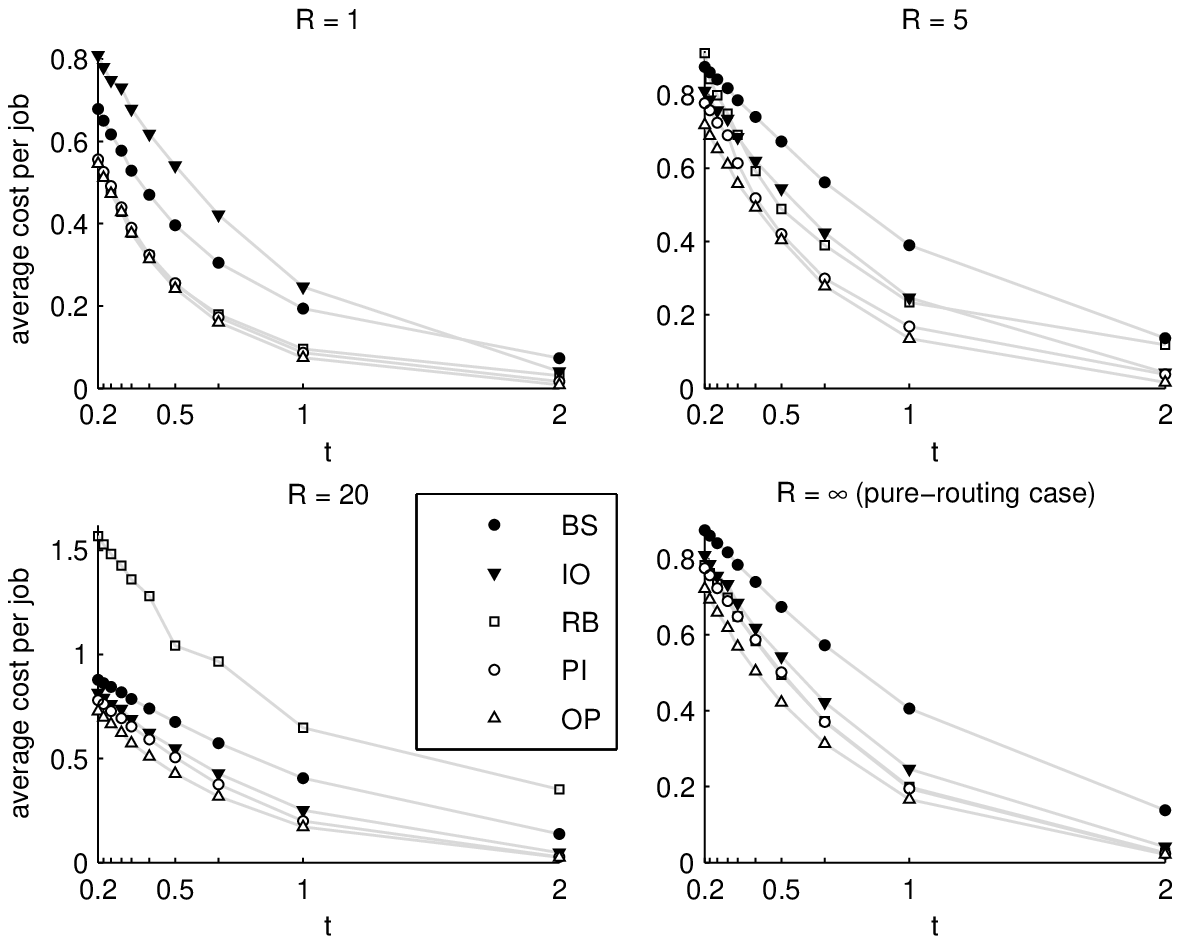}
\caption{Effect of magnitude of response-time deadline (deterministic case with $\tau \equiv t$).}
\label{fig:exp4}
\end{figure}

The fifth experiment benchmarks the performance of the policies as the response-time deadline, which is taken to be deterministic with $\tau \equiv t$, 
 ranges from a low to a high value.
We consider a base instance with server-pool sizes $\mathbf{m} = (4, 8)$, service rates $\boldsymbol{\mu} = (5, 3)$, and system load $\rho = 0.9$. 
Figure  \ref{fig:exp4} plots, for each of four values of the rejection cost $R \in \{1, 5, 20, \infty\}$, the average cost per job incurred under each policy, and under an optimal policy, against the response-time deadline $t = 1/\theta$, with  $\theta \in  \{0.1, 0.5, 1, 1.5, 2, 2.5, 3, 3.5, 4, 4.5\}$. 

For the model with admission control, the 
The first three panes in Figure \ref{fig:exp4} show that, for the model with admission control, the PI policy is consistently near optimal, with the other policies showing substantial optimality losses in some cases. Thus, the worst policy is IO for low $R$, BS for medium $R$, and RB for large $R$. For the pure-routing case, the fourth pane in the figure shows that PI and RB are both nearly optimal. 

\section{Discussion and conclusions}
\label{s:concl}
This paper has investigated the use of four policy design methods to obtain policies for control of admission and routing in a distributed soft real-time system with a parallel multicluster architecture: the static BS policy, and the dynamic IO, PI, and RB index policies.
The above analyses and numerical results reveal the following facts and insights:
\begin{enumerate}
\item All the heuristic index policies considered can be computed efficiently in linear time via first-order linear recursions.
\item The admission policy based on the IO index does not reject any job when $R \geqslant C$, even though doing so can substantially improve performance. 
\item The admission policy based on the RB index tends to grossly overreject jobs when both the system load and the rejection cost are high, since under heavy loads the RB index can take extremely high values. This causes a severe performance degradation of the RB policy in such cases. Yet, in the pure-routing case, the RB routing index yields nearly optimal policies. 
\item The optimal BS policy can be efficiently computed in the present model. Extending the range of rejection cost values $R$ beyond the case $R = 1$ considered in prior work to $R \geqslant 1$ leads to policies that may saturate one or more queues when the load is high. 
\item The optimal BS policy can be substantially outperformed by good dynamic policies.
\item Overall, the best policy among those considered for admission control and routing is the PI policy, which consistently yields nearly optimal policies. 
As for the pure-routing model, both the PI and RB index policies are consistently near optimal. 
\end{enumerate}

An interesting issue for future research would be to analyze the performance of the proposed PI index policy, or of another heuristic policy that performed as well, elucidating its asymptotic performance in instances with many servers, as the analysis of many-server models has received substantial attention in related literature (see \cite{gurvWhitt10} and references therein).

\section*{Acknowledgments}
The author acknowledges partial  support for this work by the Spanish Ministry of Science and Innovation's projects MTM2007-63140 and MTM2010-20808. He presented a preliminary outline of this work at the  Twelfth Workshop on MAthematical performance Modeling and Analysis (MAMA 2010), in New York, which appeared in the corresponding proceedings (\cite{nmmama10}).

\bibliographystyle{elsarticle-harv}



\end{document}